\documentclass{article}[11pt]
\usepackage{amsmath, amssymb, amsthm, enumerate}
\title{Multiple harmonic sums $\mathcal{H}_{\lbrace s\rbrace^{2l}=1;p-1}$ modulo $p^4$ and applications}
\author{Claire Levaillant}

\newcommand{\lb}{\lbrace}
\newcommand{\rb}{\rbrace}
\newcommand{\mpt}{\;\text{mod}\,p^3}
\newcommand{\mpd}{\;\text{mod}\,p^2}
\newcommand{\mpu}{\;\text{mod}\,p}

\newcommand{\mc}{\mathcal{C}}
\newcommand{\mo}{\mathcal{OT}}

\newcommand{\fcmcd}{\big(\mc\mb(2(p-1))\big)_1}
\newcommand{\fcmcu}{\big(\mc\mb(p-1)\big)_1}
\newcommand{\mb}{\mathcal{B}}

\newcommand{\md}{\mathcal{D}}
\newcommand{\mh}{\mathcal{H}}
\newcommand{\bn}{\binom}
\newcommand{\nts}{\negthickspace}
\newcommand{\mt}{\mathcal{T}}
\newcommand{\mS}{\mathcal{S}}
\newcommand{\as}{A^{\star}}
\newcommand{\ma}{\mathcal{A}}
\newcommand{\mcp}{\mathcal{P}}
\newcommand{\mpq}{\;\;\text{mod}\,p^4}
\newcommand{\msu}{\mathcal{S}^{'}_1}
\newcommand{\msd}{\mathcal{S}^{'}_2}
\newcommand{\mst}{\mathcal{S}^{'}_3}

\newcommand{\mbt}{\mb\mb\mb}

\newcommand{\eu}{\varphi}
\newcommand{\euq}{\varphi(p^4)}

\begin{document}
\maketitle
\begin{center}Abstract\end{center} In $1900$, Glaisher provided an expression for $(p-1)!$ modulo $p^2$ in terms of Bernoulli numbers, thus improving on Wilson's residue modulo $p$, commonly known as "Wilson's theorem". A hundred years later, Sun provides an expansion of $(p-1)!$ up to the modulus $p^3$, still in terms of Bernoulli numbers. In this paper, we investigate the modulus $p^4$ case. \\
Our method is based on studying the $p$-adic analysis of the Stirling numbers on $p$ letters with odd cycles modulo $p^4$ as well as their conjugates, the multiple harmonic sums on an even number of integers chosen amongst the first $(p-1)$ integers. By relating both entities, we also succeed to break $p$-residues of convolutions of three divided Bernoulli numbers of orders $(p-3)$ and $(p-5)$ into "smaller pieces", namely into $p$-adic residues of divided Bernoulli numbers and of convolutions of divided Bernoulli numbers. Our generalization of Wilson's theorem to the modulus $p^4$ uses in particular Miki's identity on Bernoulli numbers, Ernvall and Mets\"ankyl\"a's congruence on Bernoulli numbers and our past work on the residues of some truncated convolutions of divided Bernoulli numbers.

\section{Introduction and Main results}

Given a prime number $p$, the Stirling numbers on $p$ letters and $s$ cycles counts the number of permutations of $Sym(p)$ that decompose into a product of $s$ disjoint cycles. It can be shown that this number is also the sum of all the products of $(p-s)$ distinct integers chosen amongst the first $(p-1)$ integers. In the special case when $s=1$, we obtain the factorial $(p-1)!$.

The $p$-adic analysis of the Stirling numbers on $p$ letters goes back to the works \cite{GL1}\cite{GL} of British mathematician James Whitbread Lee Glaisher around the turn of the twentieth century. In particular, Glaisher obtains formulas for the Stirling numbers modulo $p^2$ in terms of the Bernoulli numbers when $s>1$ and generalizes Wilson's theorem to the modulus $p^2$ when $s=1$ by showing that:
$$(p-1)!=p\,B_{p-1}-p\mpd\;\;\textit{(Glaisher \cite{GL}, $1900$)}$$
A hundred years later, Chinese mathematician Zhi-Hong Sun uses Newton's formulas and his study on the generalized harmonic numbers modulo $p^3$ based on his generalizations of the Kummer congruences \cite{SU1} in order to uncover a congruence for $(p-1)!$ modulo $p^3$. Namely, he shows the following congruence.
$$(p-1)!=-\frac{pB_{p-1}}{p-1}+\frac{pB_{2(p-1)}}{2(p-1)}-\frac{1}{2}\Big(\frac{pB_{p-1}}{p-1}\Big)^2\mpt\;\;\textit{(Sun \cite{SU2}, $2000$)}$$
In \cite{SU2}, using Newton's formulas, Sun relates modulo $p^3$ the factorial $(p-1)!$ and its inverse respectively to a common convolution of divided Bernoulli numbers with ordinary Bernoulli numbers. His method does not generalize to the modulus $p^4$ as there is no such common convolution then: the perfect symmetry which exists in the modulus $p^3$ case does no longer hold in the modulus $p^4$ case. Also, going back to the modulus $p^3$, there is unfortunately no such nice general formula for all the other Stirling numbers (that is when $s\neq 1$, using our former notations). Indeed,
when applying Newton's formulas in the case when $s>1$, the Stirling numbers get related with a full convolution of divided Bernoulli numbers while their conjugates, the multiple harmonic sums on $(p-1)$ integers, rather get related to a truncated convolution of divided Bernoulli numbers. There is a way of relating both conjugate entities modulo $p^3$. This provides in turn some congruences relating the truncated convolutions to the full convolutions. The whole process involves some sophisticated techniques due to Sun which had led to the knowledge of the generalized harmonic numbers modulo $p^3$. In \cite{LEV3}, we show that there exist two other ways to bypass this knowledge. One way is namely based on using Miki's identity, an identity that is due to Japanese mathematician Hiroo Miki. His identity relates convolutions of divided Bernoulli numbers with binomial convolutions of divided Bernoulli numbers and harmonic numbers. Another way proceeds by doing some $p$-adic analysis on the polynomial $X^{p-1}+(p-1)!\in\mathbb{Z}_p[X]$, thus allowing to express the Stirling numbers modulo $p^3$ in terms of the truncated convolutions of divided Bernoulli numbers. Both ways are independent of one another.\\ The depth of the current paper relies on the fact that in the modulus $p^4$ study, there is some imbrication and complementarity between what used to be independent ways in the modulus $p^3$ case. Namely, we use Miki's identity in order to relate $(p-1)!$ and its inverse modulo $p^4$ to a common entity which, this time, is a linear combination of a convolution of three divided Bernoulli numbers whose indices sum up to $(p-1)$ and of the second residue in the $p$-adic expansion of a convolution of divided Bernoulli number also of order $(p-1)$. From there, we are able to derive $(p-1)!$ modulo $p^4$ in terms of
$$(p\,\mb_{p-1})^k,k\in\lbrace1,2\rbrace,\;p\,\mb_{2(p-1)},\;p\,\mb_{3(p-1)},\mb_{p-3}$$
This is the purpose of Theorem $3$ below.
Our method also relies on Finnish mathematicians Reijo Ernvall and Tauno Mets\"ankyl\"a's congruence which provides the difference of divided Bernoulli numbers $(\mb_{p-1+n}-\mb_n)$ modulo $p^2$ for those integers $n$ satisfying to $4\leq n\leq p-3$, in terms of the Fermat quotients in base $a$ with $1\leq a\leq p-1$.

Our work will be based on computing the Stirling numbers on $p$ letters with an odd number of cycles modulo $p^4$, and their conjugates, namely  the multiple harmonic sums with an even number of integers, each carrying the exponent one.
As a matter of fact, we will use these computations only with some specific values for the odd number of cycles in the first case or for the even number of integers in the latter case. But we work in the more general setting in order to allow for future developments in forthcoming papers.
By working out some special cases of $s>1$ such as $s=3$ resp $s=5$, we are able to reduce modulo $p$ a convolution of three divided Bernoulli numbers of order $(p-3)$, resp $(p-5)$ to a convolution of two divided Bernoulli numbers of the same order, as in Theorem $1$ below.\\
In \cite{SU2}, Sun provides without proof the generalized harmonic numbers $H_{p-1,k}$ modulo $p^4$, but only for those $k$ with $1\leq k\leq p-5$. In the current paper, we need and thus treat the cases $k\in\lb p-3,p-1\rb$ as well. This is the purpose of Theorem $2$ below. In order to do so, we adapt Sun's proof done for the modulus $p^3$ to the modulus $p^4$. Sun's proof uses Euler's totient function, sums of powers of integers and his Kummer type congruences of \cite{SU0} which he later generalizes in \cite{SU2}.

Throughout the paper we denote the divided Bernoulli numbers by $\mb$ and the ordinary Bernoulli numbers by $B$.
By definition,
$$\mb_i:=\frac{B_i}{i}\;\forall\,i\geq 1$$
We will denote the $k$-th harmonic number by $\mathcal{H}_k$ and the generalized harmonic numbers by $H_{p-1,k}$ or simply $H_k$. We let $\varphi$ denote Euler's totient function. For an integer $n$, $\eu(n)$ is the order of the group of units of $\mathbb{Z}/n\mathbb{Z}$. We recall that $\eu$ is multiplicative and $\eu(p^{\alpha})=p^{\alpha-1}(p-1)$ when $p$ is prime. When $x$ is a $p$-adic integer, we will denote by $(x)_i$ the $i$-th coordinate of its Hensel's expansion and we shall also use the same notation for the nonnegative powers of $p$ coordinates when $x$ belongs to $\mathbb{Q}_p$ instead of $\mathbb{Z}_p$. \textbf{However, we make an exception for $(pB_{p-1})_1$ which will exceptionally denote the Agoh-Giuga quotient, that is:} $$\mathbf{pB_{p-1}=-1+p(pB_{p-1})_1}$$ We will denote by $\mc\mb(n)$ the convolution of divided Bernoulli number of order $n$, with starting index $2$ and ending index $n-2$. And so by definition,
$$\mc\mb(n):=\sum_{i=2}^{n-2}\mb_i\mb_{n-i}$$
We will denote by $\mbt(n)$ the sum of products of three divided Bernoulli numbers whose indices sum up to $n$. Again, it is understood that all the indices are greater than or equal to $2$.\\

\textbf{Unless otherwise specified, throughout the paper $p$ is a given prime with $p>7$.} \\

We are now ready to state the main results. \\

\noindent The first main theorem deals with convolutions of three divided Bernoulli numbers of respective orders $p-3$ and $p-5$.

\newtheorem{Theorem}{Theorem}

\begin{Theorem} Let $w_p$ denote the Wilson quotient, and so we have by \cite{LEV1}:
\begin{eqnarray*}(w_p)_0&=&\left(\begin{array}{l}\frac{p\,B_{p-1}+1}{p}\end{array}\right)_0-1\;\;\mpu\;\mathbb{Z}_p\\
(w_p)_1&=&\bigg(1+\frac{p}{2}-\frac{3}{2}p^2+(2p+1)pB_{p-1}-\frac{1}{2}pB_{2(p-1)}-\frac{1}{2}p^2B_{p-1}^2\bigg)_2
\end{eqnarray*}
%(w_p)_1&=&\begin{split}&\left(\begin{array}{l}\left(\begin{array}{l}\frac{p\,B_{p-1}+1}{p}\end{array}\right)_0-1\end{array}\right)_1
%+\left(\begin{array}{l}\frac{p\,B_{p-1}+1}{p}\end{array}\right)_1
%-\frac{1}{2}\left(\begin{array}{l}\frac{p\,B_{2(p-1)}+1}{p}\end{array}\right)_1\\&+2\left(\begin{array}{l}\frac{p\,B_{p-1}+1}{p}\end{array}\right)_0
%-\frac{1}{2}
%\left(\begin{array}{l}\frac{p\,B_{p-1}+1}{p}\end{array}\right)_0^2-\frac{3}{2}\,\qquad\qquad\qquad\mpu\mathbb{Z}_p\end{split}\end{eqnarray*}
%-\frac{3}{2}+2(p\,B_{p-1})_1+(p\,B_{p-1})_2-\frac{1}{2}\big(p^2\,B_{p-1}^2+p\,B_{2(p-1)}\big)_2\,\mpu

The following congruences hold.\\\\
$(i)$ $\mbt(p-3)$ is congruent modulo $p$ to
$$\begin{array}{l}
3\Big(\mc\mb(p-3)\Big)_1-\frac{11}{4}\mb_{p-5}+6(\mb_{3p-5}-3\mb_{2p-4}+2\mb_{p-3})_2\\
+\Big[6(\mb_{p-3}-\mb_{2p-4})_1\Big]_1+\Big[3\big(2(\mb_{2p-4}-\mb_{p-3})_1+2(w_p)_0\mb_{p-3}\big)_0\Big]_1\\
+(w_p)_1\Big[-6\mb_{p-3}\Big]_0+(w_p)_0\Big[6(\mb_{2p-4}-2\mb_{p-3})\Big]_1+\Big[(w_p)_0(-6\mb_{p-3})_0\Big]_1\\
+\Big\lbrace\Big[6(\mb_{p-3}-\mb_{2p-4})_1\Big]_0+\Big[-6(w_p)_0\mb_{p-3}\Big]_0\\
+\Big[3\big(2(\mb_{2p-4}-\mb_{p-3})_1+2(w_p)_0\mb_{p-3}\big)_0\Big]_0\Big\rbrace_1
\end{array}$$

%(i)\;\mbt(p-3)=&\;\\%\Big(3\mc\mb(p-3)\Big)_1-\frac{11}{4}(\mb_{p-5})_0+\Big[6(\mb_{3p-5}-3\mb_{2p-4}+2\mb_{p-3})\Big]_2\\
%&+(w_p)_1(-6\mb_{p-3})_0+(w_p)_0\Big[6(\mb_{2p-4}-2\mb_{p-3})\Big]_1+\Big[(w_p)_0(-6\mb_{p-3})_0\Big]_1\\
%&+\Big\lbrace(\Big[6(\mb_{2p-4]-\mb_{p-3})_1+6(w_p)_0\mb_{p-3}\Big]_0+\Big[6(\mb_{3p-5}-3\mb_{2p-4}+2\mb_{p-3})\Big]_1+\Big[-6(w_p)_0\mb_{p-3}\Big]_0 %\Big\rbrace_1\\&
%\qquad\qquad\qquad\qquad\qquad\qquad\mpu
%\end{split}\end{equation*}
\noindent $(ii)$ $\mbt(p-5)$ is congruent modulo $p$ to
$$\begin{array}{l}
3\big(\mc\mb(p-5)\big)_1-\frac{15}{4}\mb_{p-7}+\big[3(\mb_{p-3})_0^2]_1+6\mb_{p-3}\big(2\mb_{p-3}-\mb_{2p-4}\big)_1\\
+6\big(\mb_{3p-7}-3\mb_{2p-6}+2\mb_{p-5}\big)_2+\Big[6(\mb_{p-5}-\mb_{2p-6})_1\Big]_1\\
+(w_p)_0(-6\mb_{p-5})_1+(w_p)_1(-6\mb_{p-5})_0\\
+\Big((w_p)_0(-6\mb_{p-5})_0\Big)_1+\Big(-6(w_p)_0^2\mb_{p-5}\Big)_0\\
+\Big[3(w_p)_0(2(w_p)_0\,\mb_{p-5}+2(\mb_{2p-6}-\mb_{p-5})_1-\mb_{p-3}^2)_0\Big]_0\\
+\Big[3\big(2(w_p)_0\,\mb_{p-5}+2(\mb_{2p-6}-\mb_{p-5})_1-\mb_{p-3}^2\big)_0\Big]_1\\
+\Big\lbrace\Big[3(\mb_{p-3})_0^2\Big]_0+\Big[6(\mb_{p-5}-\mb_{2p-6})_1\Big]_0+\Big[-6(w_p)_0(\mb_{p-5})_0\Big]_0\\
+\Big[3\big(2(w_p)_0\,\mb_{p-5}+2(\mb_{2p-6}-\mb_{p-5})_1-\mb_{p-3}^2\big)_0\Big]_0\Big\rbrace_1
%&\Big[-\big((w_p)_1(\mb_{p-5})_0+(w_p)_0(\mb_{p-5})_1\big)\Big)\\
\end{array}$$
\end{Theorem}

\noindent The next main theorem provides the expansion of the generalized harmonic numbers $H_{p-1,p-3}$ and $H_{p-1,p-1}$ up to the modulus $p^4$. These are respectively needed for the computations of the multiple harmonic sum on $p-1$ integers of order $p-3$, resp. $p-1$, modulo $p^4$.

\begin{Theorem} Modulo $p^4$ the following congruences hold.
\begin{equation*}\begin{split}
%H_{p-3}=-\frac{1}{2}\big(p-3+\frac{11}{2}p^3\big)B_{2p}+p(3-p)\big(-3\,\mb_{p+1}-\frac{p-1}{4}\big)+p^2(2(p-1)-pB_{p-1})+\frac{5}{6}p^3\\
H_{p-3}\equiv\frac{37p^3}{12}-3p^2+\frac{3p}{4}-p^2(pB_{p-1})+3p&(p-3)\mb_{p+1}\\&-\frac{1}{2}\Big(p-3+\frac{11p^3}{2}\Big)B_{2p}\end{split}\end{equation*}
\begin{equation*}\begin{split}H_{p-1}\equiv\frac{13p^3+12p-12}{3}+\frac{19p^3-12}{2}\,pB_{p-1}+&(4-7p^3)\,pB_{2(p-1)}\\&+\frac{11p^3-6}{6}\,pB_{3(p-1)}
\end{split}\end{equation*}
\end{Theorem}

\noindent In \cite{SU2}, Sun provides without proof the reduction modulo $p^4$ of the generalized harmonic numbers $H_{p-1,k}$, but the cases $k=p-3$ and $k=p-1$ of interest here are omitted.\\

The last main theorem offers a generalization of Wilson's theorem to the modulus $p^4$, thus pushing the successive p-adic expansions of Wilson, Glaisher and Sun one p power further.

\begin{Theorem} Let $p\geq 13$ be a prime. P-adic expansion of $(p-1)!$ up to the modulus $p^4$. \\In what follows, the other term and residue $\widehat{\mo}$ gets provided in Appendix A.
\begin{equation*}\begin{split}
(p&-1)!=\frac{1}{3}\Bigg\lbrace -p^3\mb_{p-3}+\frac{p^2}{2}(2p\mb_{2(p-1)}-p^2\mb_{p-1}^2)_2(3-p\big(1+(pB_{p-1})_1\big)\\
&+4-10p\mb_{p-1}+8p\mb_{2(p-1)}-3p\mb_{3(p-1)}-\frac{p^3}{2}\widehat{\mo}\\
%&-\Big(p^3\mb_{p-1}^3+p^2\mb_{p-1}^2-2(p\mb_{2(p-1)})(p\mb_{p-1})-1\Big)\Big(3(p\mb_{p-1}-1)-p\mb_{2(p-1)}-\frac{1}{2}p^2\mb_{p-1}^2\Big)\\
&+\Big(p\mb_{p-1}\big[2p\mb_{2(p-1)}-p^2\mb_{p-1}^2\mpt\big]-p^2\mb_{p-1}^2+1\Big)\Big(\frac{1}{(p-1)!}\mpt\Big)\Bigg\rbrace\\
%&-p^3\mb_{p-3}-\frac{p^3}{6}(\mb_{2p-4}-\mb_{p-3})_1+\frac{p^3}{30}(\mb_{2p-6}-\mb_{p-5})_1\Bigg\rbrace\\
&\qquad\qquad\qquad\qquad\qquad\qquad\qquad\qquad\qquad\qquad\qquad\qquad\qquad\;\;\;\;\mpq
\end{split}\end{equation*}
And a result of \cite{LEV3} states that:
$$\frac{1}{(p-1)!}=3(p\mb_{p-1}-1)-p\mb_{2(p-1)}-\frac{1}{2}p^2\mb_{p-1}^2\mpt$$

%&(p-1)!=\frac{1}{3}\Big\lbrace 1-p^3\Big(2+(\mb_{p-1})_0(\mb_{p-1})_1+(\mb_{p-1})_2\Big)\\
%&+\frac{1}{2}p^5\mb_{p-1}^5-\frac{5}{2}p^4\mb_{p-1}^4+\frac{p-1}{2}p^3\mb_{p-1}^3+\frac{p+3}{2}p^2\mb_{p-1}^2-7p\mb_{p-1}\\
%&+\Big(7p^2\mb_{p-1}^2-(p+5)p\mb_{p-1}-p+9\Big)p\mb_{2(p-1)}-2(p\mb_{p-1})\,(p^2\mb_{2(p-1)}^2)-3p\mb_{3(p-1)}\\
%&-p^3\Big(2(p^2\mb_{p-1}^2)_0+(p^2\mb_{p-1}^2)_1+\frac{1}{2}(p^2\mb_{p-1}^2)_2\Big)-p^3\big(\mh_{2(p-1)}\mb_{2(p-1)}-2\,\mh_{p-1}\mb_{p-1}\big)_1\\&
%-\frac{p^2}{6}\mb_{2(p-2))}+\frac{p^2}{30}\mb_{2(p-3)}+\bigg(\frac{p^2}{6}+p^3\Big((\mb_2)_0+(\mb_2)_1-\frac{2}{3}\Big)\bigg)\,\mb_{p-3}-\bigg(\frac{p^2}{10}+\frac{p^3}{15}\Big((p\,B_{p-1})_1+\frac{1}{4}\Big)\bigg)\mb_{p-5}\Bigg\rbrace\\
%&\qquad\qquad\qquad\qquad\qquad\qquad\qquad\qquad\qquad\qquad\qquad\qquad\qquad\qquad\qquad\qquad\qquad\qquad\qquad\;\mpq\end{split}\end{equation*}
\end{Theorem}

The paper is organized as follows. \\
In the next section $\S\,2$, we will recall some past results of importance here.

In $\S\,3$ we will compute the Stirling numbers on $p$ letters and with an odd number $s$ of cycles to the modulus $p^4$. Moreover, we will relate modulo $p^4$ these numbers and the multiple harmonic sums on an even number of integers chosen amongst the first $(p-1)$ integers in a certain way, thus generalizing an idea of \cite{LEV3} to the next $p$-power. We will then use these congruences in order to compute the convolutions from Theorem $1$ modulo $p$ in the cases when $s=3$ and $s=5$. Along the way we will need to compute the multiple harmonic sums on respectively two and four integers modulo $p^4$.

Finally in $\S\,4$, we study the multiple harmonic sums on an even number of integers chosen amongst the first $(p-1)$ integers modulo $p^4$ without any reference to the Stirling numbers. From this second congruence modulo $p^4$ on the multiple harmonic sums,
we show how to derive $(p-1)!$ mod $p^4$ in the case when $s$ is chosen to be $1$. This is the part requesting a fair amount of work.\\

\section{Some useful preliminary results}

An efficient way in order to compute the Stirling numbers or the multiple harmonic sums modulo powers of primes was developed in \cite{SU2} and relies on the use of Newton's formulas.
Already Glaisher around the turn of the twentieth century was using Newton's formulas in order to related the generalized harmonic numbers and the Stirling numbers modulo the small powers $p^2$ or $p^3$.\\

We recall these old and useful formulas below.

\newtheorem{Result}{Result}
\begin{Result}(Newton's formula, see for instance \cite{JAC}). \\
Let $x_1,\dots,x_m$ be complex numbers.\\
Let $\mcp_k=x_1^k+\dots +x_m^k$ and $\ma_k=\sum_{1\leq i_1<\dots<i_k\leq m}x_{i_1}\dots x_{i_k}$. \\\\
Then, for $k=0,1,\dots,m$ we have:
\begin{equation*}\mcp_k-\ma_1\mcp_{k-1}+\ma_2\mcp_{k-2}+\dots+(-1)^{k-1}\ma_{k-1}\mcp_1+(-1)^{k}k\ma_k=0\end{equation*}
\end{Result}

\noindent Sums of particular importance throughout the paper will be when the $x_i$'s above are chosen to be integers amongst the first $(p-1)$ integers, or reciprocals of those. In the first situation, the $\ma_k$'s and $\mcp_k$'s are respectively denoted by $A_k$ and $S_k$ and respectively called unsigned Stirling numbers of the first kind and sums of powers. In the second situation, the $\ma_k$'s and $\mcp_k$'s are respectively denoted by $\as_k$ and $H_k$ and respectively called multiple harmonic sums and generalized harmonic numbers. In a more general context the generalized harmonic numbers $H_k$ can also be denoted by $H_{p-1,k}$ in order to signify their order. \\

\noindent The $p$-adic expansions of $A_k$, $\as_k$, $S_k$ and $H_k$ can all be expressed in terms of Bernoulli numbers. P-adically, the well-known theorem on Bernoulli numbers that is due to Von Staudt and independently Clausen is of crucial importance. We state below a quite useful consequence of it.

\begin{Result}(Consequence of Von Staudt-Clausen's theorem, \cite{CL} and independently \cite{VS} $1840$).\\
The denominator of the $i$th Bernoulli number $B_i$ consists of a product of primes $p$ all with multiplicities $1$, such that $p-1$ divides $i$.
\end{Result}

\noindent By using the Bernoulli formula for the sums of powers, it is possible to derive the following congruence. This is for instance achieved in \cite{SU2} or \cite{LEV1}.
\begin{Result} (Sun \cite{SU2}, $2000$) Let $k$ be an integer with $k\geq 2$. Then, we have:
$$S_k=p\,B_k+\frac{p^2}{2}\,k\,B_{k-1}+\frac{p^3}{6}\,k(k-1)\,B_{k-2}+\frac{k(k-1)(k-2)}{24}p^4\,B_{k-3}\;\;\text{mod}\,p^4$$
\end{Result}

\noindent The oldest congruences on harmonic numbers and generalized harmonic numbers go back to Wolstenholme. The result below is known as "Wolstenholme's theorem" although Wolstenholme's actual theorem dealt with binomial coefficients and the result below only stood as an intermediate result.
\begin{Result} (due to Wolstenholme, \cite{WO} $1862$).
\begin{eqnarray*}H_1=H_{p-1,1}&=&0\mpd\\
H_2=H_{p-1,2}&=&0\mpu\end{eqnarray*}
\end{Result}

\noindent The result got more recently generalized by Bayat in \cite{BA}, where he also deals with the other generalized harmonic numbers $H_{p-1,k}$ for $3\leq k\leq p-3$. However, the cases $k=p-2$ and $k=p-1$ are left unstudied. By generalizing Kummer's congruences (see Result $12$ below) to the modulus $p^2$ and using other congruences proven by himself in \cite{SU1}, Sun finds the $p$-adic expansion of all the generalized harmonic numbers $H_{p-1,k},\,1\leq k\leq p-1$ up to the modulus $p^3$ and even to the modulus $p^4$ in the case when $1\leq k\leq p-5$, as he is able to generalize the Kummer's congruences one $p$ power even further. In the current paper, we deal with the extra cases $k=p-3$ and $k=p-1$, see Theorem $2$ of $\S\,1$.\\
We first state the reduction of the modulus $p^3$ to the modulus $p^2$, then state the moduli $p^3$ and $p^4$ cases.

\begin{Result} (due to Glaisher, \cite{GL} $1900$). \\
Let $p$ be a prime greater than $3$ and $k\in\lb 1,2,\dots,p-1\rb$. Then, \\
$$\sum_{a=1}^{p-1}\frac{1}{a^k}=\begin{cases}\frac{k}{k+1}\,p\,B_{p-1-k}\qquad\qquad(mod\,p^2)&\text{if $k<p-1$}\\
-p\,B_{p-1}+2(p-1)\;\;\;\;\;\,(mod\,p^2)&\text{if $k=p-1$}
\end{cases}$$
\end{Result}

\begin{Result} (due to Sun, \cite{SU2} $2000$). \\Let $p$ be a prime greater than $3$. Then,\\
(a) If $k\in\lbrace 1,2,\dots,p-4\rbrace$ then
$$\sum_{a=1}^{p-1}\frac{1}{a^k}=\begin{cases}\frac{k(k+1)}{2}\frac{B_{p-2-k}}{p-2-k}p^2\qquad\qquad\;\mpt&\text{if $k$ is odd}\\
k\bigg(\frac{B_{2p-2-k}}{2p-2-k}-2\frac{B_{p-1-k}}{p-1-k}\bigg)p\;\mpt&\text{if $k$ is even}
\end{cases}$$
(b) $$\sum_{a=1}^{p-1}\frac{1}{a^{p-3}}=\bigg(\frac{1}{2}-3B_{p+1}\bigg)p-\frac{4}{3}p^2\mpt$$
(c) $$\sum_{a=1}^{p-1}\frac{1}{a^{p-2}}=-\big(2+p\,B_{p-1}\big)p+\frac{5}{2}p^2\mpt$$
(d) $$\sum_{a=1}^{p-1}\frac{1}{a^{p-1}}=p\,B_{2p-2}-3p\,B_{p-1}+3(p-1)\mpt$$
%In our paper, we will name these numbers generalized harmonic numbers and we will denote them by $H_{p-1,k}$ or simply $H_k$, using standard notations.
\end{Result}

\begin{Result} (due to Sun, \cite{SU2} $2000$). \\
Let $p>5$ be a prime and $k\in\lb 1,2,\dots,p-5\rb$. Then, modulo $p^4$, we have:\\
$$\sum_{a=1}^{p-1}\frac{1}{a^k}\equiv\begin{cases} -k(\mb_{3p-3-k}-3\mb_{2p-2-k}+3\mb_{p-1-k})p-\binom{k+2}{3}p^3\mb_{p-3-k}&\text{if $2|k$}\\
-\binom{k+1}{2}(\mb_{2p-3-k}-2\,\mb_{p-2-k})\,p^2&\text{if $2\not|k$}\end{cases}$$
\end{Result}
\noindent The next results deal with the Stirling numbers modulo $p^3$ and the multiple harmonic sums modulo $p^3$. They are respectively stated in \cite{LEV1} and \cite{LEV3}.
\begin{Result}(Stirling numbers modulo $p^3$, cf \cite{LEV1}, preprint $2019$). \\Let $k$ be an integer with $1\leq k\leq \frac{p-1}{2}$. We have,
\begin{eqnarray*}
(k\leq\frac{p-3}{2})\;A_{2k+1}\negthickspace\negthickspace&=&\negthickspace\!\!\frac{p^2}{2}\frac{2k+1}{2k}\,B_{2k}\qquad\qquad\qquad\qquad\qquad\!\text{mod}\,p^3\\
&&\notag\\
(k\neq 1)\;A_{2k}\negthickspace\negthickspace&=&\negthickspace\!\!-\frac{1}{2k}\Bigg(p\,B_{2k}-p^2\,\sum_{r=1}^{ k-1}\frac{B_{2r}B_{2k-2r}}{2r}\Bigg)\;\;\;\text{mod}\,p^3\\
&&\notag\\
\text{and}\;A_1\negthickspace\negthickspace&=&\negthickspace\!\!\frac{p(p-1)}{2}\qquad\qquad\qquad\qquad\qquad\qquad\text{mod}\,p^3\\
&&\notag\\
(p\geq 5)\;\text{and}\;A_2\negthickspace\negthickspace&=&\negthickspace\!\!\frac{1}{2}\bigg(-\frac{p}{6}+\frac{3\,p^2}{4}\bigg)\qquad\qquad\qquad\qquad\;\text{mod}\,p^3
\end{eqnarray*}
\end{Result}

\begin{Result} (Multiple harmonic sums modulo $p^3$, \cite{LEV3} $2020$). \\
Let $\mt\mc\mb(p+1-k,p-3)$ denote the following truncated convolution of divided Bernoulli numbers:
$$\mt\mc\mb(p+1-k,p-3)=\sum_{i=p+1-k}^{p-3}\mb_i\mb_{2(p-1)-k-i}$$
$$\as_k\equiv\begin{cases}\frac{k+1}{2}\,p^2\,\mb_{p-2-k}&\text{if $k$ is odd and $k\leq p-4$}\\
&\\
\frac{p}{2}-p^2+\frac{1}{2}(p\,B_{p-1})_1p^2&\text{if $k=p-2$}\\
&\\
\,p\,\big(2\,\mb_{p-1-k}-\mb_{2(p-1)-k}\big)&\\&\\\qquad\qquad\;\;+\frac{p^2}{2}\mt\mc\mb(p+1-k,p-3)&\text{if $k$ is even and $k\leq p-5$}\\
&\\
\frac{p}{12}-\frac{11p^2}{24}+\frac{p^2}{12}(p\,B_{p-1})_1&\text{if $k=p-3$}\\
&\\
3(p\,\mb_{p-1}-1)-p\,\mb_{2(p-1)}-\frac{1}{2}p^2\,\mb_{p-1}^2&\text{if $k=p-1$}
\end{cases}$$
\end{Result}
\noindent The following congruences relate some full convolutions of divided Bernoulli numbers to some truncated convolutions of divided Bernoulli numbers. They will be useful in the present paper. In \cite{LEV3} we present three independent ways to derive them, like discussed in the introduction.
\begin{Result} (Truncated convolutions in terms of full convolutions modulo $p^3$, \cite{LEV3} $2020$).\\
(i) Assume that $4\leq 2n\leq p-7$.
\begin{equation*}
\begin{split}
\frac{p^2}{2}\sum_{i=p+1-2n}^{p-3}\mb_i\mb_{2(p-1)-2n-i}=&-\frac{p^2}{2}\sum_{i=2}^{p-1-2n-2}\mb_i\mb_{p-1-2n-i}\\
&+p\,\big(\mb_{2(p-1)-2n}-\mb_{p-1-2n}\big)\\
&+p^2\bigg((p\,B_{p-1})_1-1\bigg)\,\mb_{p-1-2n}\mpt
\end{split}
\end{equation*}
(ii) Case $2n=p-3$.
\begin{equation*}
p^2\sum_{i=4}^{p-3}\mb_i\mb_{p+1-i}=2p\,\mb_{p+1}+p^2\mb_2+2p\mb_2(p\,B_{p-1})\mpt
\end{equation*}
(iii) Case $2n=p-5$.
\begin{equation*}
p^2\sum_{i=6}^{p-3}\mb_i\mb_{p+3-i}=\frac{7}{720}\,p^2+2p\,\mb_{p+3}+2p\,\mb_4(p\,B_{p-1})\mpt
%\sum_{i=6}^{p-3}\mb_i\mb_{p+3-i}=\Bigg(\bigg(2\mh_3-\frac{3}{2}\bigg)\mb_4+B_2^2\Bigg)p^2+2p\,\mb_{p+3}+2p\mb_4(p\,B_{p-1})\mpt
\end{equation*}
\end{Result}
\noindent In \cite{LEV3}, it is also shown that full convolutions of orders $p-1$ and $p-3$ and $p-5$ can be reduced modulo $p$ as follows.
\begin{Result} (Full convolutions of divided Bernoulli numbers of respective orders $p-1$, $p-3$ and $p-5$ modulo $p$, \cite{LEV3} $2020$).\\
\begin{eqnarray*}
p^2\,\mc\mb(p-1)&=&2\,p\,\mb_{2(p-1)}-p^2\mb_{p-1}^2\qquad\qquad\qquad\qquad\;\mpt\\
\mc\mb(p-3)&=&2\,w_p\,\mb_{p-3}+2\,(\mb_{2p-4}-\mb_{p-3})_1\qquad\qquad\mpu\\
\mc\mb(p-5)&=&2\,w_p\,\mb_{p-5}+2\,(\mb_{2(p-3)}-\mb_{p-5})_1-\mb_{p-3}^2\,\!\mpu
\end{eqnarray*}
\end{Result}

The following three results state congruences concerning Bernoulli numbers $B_{k(p-1)+b}$, all of which got mentioned earlier.

\begin{Result}(Kummer's congruences, $1850$ \cite{KU}, followed by Sun's generalization to the moduli $p^2$ and $p^3$, $2000$ \cite{SU2}).\\
Let $p$ be an odd prime and let $b>0$ be an even integer such that $p-1\not|\,b$. Then, we have for all nonnegative integer $k$,
$$\frac{B_{k(p-1)+b}}{k(p-1)+b}\equiv\frac{B_b}{b}\;\;\text{mod}\,p$$
$$\frac{B_{k(p-1)+b}}{k(p-1)+b}\equiv k\frac{B_{p-1+b}}{p-1+b}-(k-1)(1-p^{b-1})\frac{B_b}{b}\;\;\text{mod}\,p^2$$
$$\frac{B_{k(p-1)+b}}{k(p-1)+b}\equiv \binom{k}{2}\frac{B_{2(p-1)+b}}{2(p-1)+b}-k(k-2)\frac{B_{p-1+b}}{p-1+b}+\binom{k-1}{2}(1-p^{b-1})\frac{B_b}{b}$$
$$\;\;\qquad\qquad\qquad\qquad\qquad\qquad\qquad\qquad\qquad\qquad\qquad\qquad\qquad\qquad\qquad\text{mod}\,p^3$$
We will denote these congruences respectively by $(K1)$, $(K2)$ and $(K3)$.
\end{Result}

\noindent When $k=1$, Sun's generalization in the modulus $p^2$ is a trivial congruence. Likewise, when $k=2$, Sun's generalization in the modulus $p^3$ case is trivial. In the first case however, some knowledge gets provided by Ernvall and Mets\"ankyl\"a. While studying the irregular primes, Ernvall and Mets\"ankyl\"a have namely shown the following congruence, where $q_a$ denotes the Fermat quotient in base $a$.
\begin{Result} (Ernvall and Mets\"ankyl\"a \cite{EM}, $1991$)
Let $n$ be an even integer with $4\leq n\leq p-3$. Then, we have:
\begin{equation*}
\mb_{p-1+n}=\mb_{n}-\frac{p}{2}\,\sum_{a=1}^{p-1}q_a^2\,a^n\;\;\mpd
\end{equation*}
\end{Result}
\noindent Their congruences will be useful in the present work.

\noindent Other Kummer type congruences of importance were developed in \cite{SU1}. We will only use some special cases which we state below.
\begin{Result} (special cases of Corollary $4.2$ of Sun \cite{SU1} $1997$).\\
Applied with $b=0$, $n=2,3$, $p>3$ and $k\geq 1$.
\begin{eqnarray*}
p\,B_{k(p-1)}&=&-(k-1)(p-1)+k\,p\,B_{p-1}\;\mpd\\
p\,B_{k(p-1)}&=&\frac{(k-2)(k-1)}{2}(p-1)-k(k-2)\,pB_{p-1}+\frac{k(k-1)}{2}\,p\,B_{2p-2}\\
&&\qquad\qquad\qquad\qquad\qquad\qquad\qquad\qquad\qquad\qquad\qquad\;\;\mpt
\end{eqnarray*}
Applied with $b=0$, $n=4$, $p>5$ and $k\geq 1$.
\begin{equation*}\begin{split}
p\,B_{k(p-1)}=&-\frac{(k-1)(k-2)(k-3)}{6}(p-1)+\frac{(k-2)(k-3)}{2}k\,pB_{p-1}\\&-\frac{k(k-1)(k-3)}{2}\,pB_{2(p-1)}
+\frac{k(k-1)(k-2)}{6}\,pB_{3(p-1)}\mpq
\end{split}\end{equation*}
We will denote these congruences after Sun by $(S2)$, $(S3)$ and $(S4)$ respectively.
\end{Result}
Last, we recall from \cite{LEV3} the following congruence.
\begin{Result} (Inverse of the factorial modulo $p^3$, \cite{LEV3} $2020$)
\begin{eqnarray*}
%(p-1)!&=&-p\,\mb_{p-1}+p\,\mb_{2(p-1)}-\frac{1}{2}p^2\,\mb_{p-1}^2\qquad\;\;\mpt\\
\frac{1}{(p-1)!}&=&3(p\,\mb_{p-1}-1)-p\,\mb_{2(p-1)}-\frac{1}{2}\,p^2\,\mb_{p-1}^2\;\mpt
\end{eqnarray*}
\end{Result}

\section{Reduction of $\mbt(p-3)$ and $\mbt(p-5)$ mod $p$}

This part is based on Newton's formulas applied to the Stirling numbers and to the sums of powers. \\Let $k\in\lbrace 1,\dots,p-1\rbrace$. By applying Result $1$ with $x_i=i,\;1\leq i\leq p-1$, we obtain:
\begin{equation}
A_k=\frac{(-1)^{k-1}}{k}\Big(S_k+\sum_{r=1}^{k-1}(-1)^r\,A_r\,S_{k-r}\Big)
\end{equation}
Assume $k$ is even. If $r$ is odd, then $k-r$ is odd and so $p^2$ divides $S_{k-r}$, unless $r=k-1$. Also we know from Result $8$ that $p^2$ divides $A_3,A_5,\dots,A_{p-2}$.
Then, modulo $p^4$, only the even indices participate into the sum of $(1)$, as well as the indices $r=1$ and $r=k-1$. Moreover, if $r$ is even,
\begin{equation}
A_r=\begin{cases}\frac{p^2}{2}\mc\mb(r)-p\,\mb_{r}\;\;\mpt&\text{if $r\neq 2$}\\
&\\
\frac{-p}{12}+\frac{3p^2}{8}\qquad\;\;\;\;\;\!\mpt&\text{if $r=2$}\end{cases}
\end{equation}
So, applying $(1)$ with $k=2n\geq 4$ yields:
\begin{equation}\begin{split}
A_{2n}=-\frac{1}{2n}\Bigg(&p\,B_{2n}+\frac{p^3}{3}n(2n-1)\,B_{2n-2}-A_1S_{2n-1}-A_{2n-1}S_1\\&
+\bigg(-\frac{p}{12}+\frac{3p^2}{8}\bigg)S_{2n-2}+\bigg(\frac{p}{6}-\frac{p^2}{2}\bigg)A_{2n-2}\\&+\sum_{r=4}^{2n-4}p\,B_{2n-r}\,\bigg
(\frac{p^2}{2}\mc\mb(r)-p\,\mb_r\bigg)\Bigg)\qquad\mpq
\end{split}\end{equation}
Further, we have, \begin{eqnarray}A_{2n-1}S_1&=&\frac{p^3(p-1)}{4}\frac{2n-1}{2n-2}\,B_{2n-2}\;\;\;\mpq\qquad\text{(cf Result $8$)} \\
A_1S_{2n-1}&=&\frac{p^3(p-1)}{4}(2n-1)\,B_{2n-2}\;\mpq\qquad\text{(cf Result $3$)}\end{eqnarray}
Furthermore, imposing $2n>4$, we have:
\begin{eqnarray}
S_{2n-2}&=&p\,B_{2n-2}\;\qquad\qquad\qquad\;\;\;\;\mpt\\
A_{2n-2}&=&\frac{p^2}{2}\mc\mb(2n-2)-p\,\mb_{2n-2}\;\mpt
\end{eqnarray}
We gather the final form inside the following statement, where we also treated the cases $2n=2$ and $2n=4$.
\newtheorem{Proposition}{Proposition}
\begin{Proposition} Stirling numbers $c(p,p-2n),\;2\leq 2n\leq p-1$ modulo $p^4$.\\\\
(i) Assume $6\leq 2n\leq p-1$. Then,
\begin{equation}\begin{split}A_{2n}=
-\frac{1}{2n}\Bigg(p\,&B_{2n}+\frac{n(16n^2-12n+5)}{24(n-1)}\,p^3\,B_{2n-2}\\
&-p^2\,n\mc\mb(2n)+\frac{p^3}{12}\,\mc\mb(2n-2)\\&+\frac{p^3}{2}\sum_{r=4}^{2n-4}B_{2n-r}\mc\mb(r)\Bigg)\;\;\mpq\end{split}\end{equation}
(ii) (Case $2n=4$).$$A_4=\frac{1}{8}\Big(\frac{p}{15}+\frac{p^2}{36}-\frac{5p^3}{4}\Big)\;\;\mpq$$
(iii) (Case $2n=2$).
$$A_2=\frac{p(1-2p)(5p-2)}{24}\;\;\mpq$$
\end{Proposition}
Let's investigate the terms in $(i)$ more closely. First and foremost, we have:
$$\sum_{r=4}^{2n-4}B_{2n-r}\mc\mb(r)=\sum_{r=2}^{2n-4}B_r\mc\mb(2n-r)-\frac{1}{6}\mc\mb(2n-2)$$
It follows that
\begin{equation}\begin{split}
A_{2n}=-\frac{1}{2n}\Bigg(&p\,B_{2n}+\frac{n(16n^2-12n+5)}{24(n-1)}\,p^3\,B_{2n-2}\\&
-p^2\,n\mc\mb(2n)+\frac{p^3}{2}\sum_{r=2}^{2n-4}B_r\mc\mb(2n-r)\Big)\mpq\end{split}
\end{equation}
Further, by considering thrice the sum in the sum of the second row of $(9)$, we obtain:
$$\sum_{r=2}^{2n-4}B_r\mc\mb(2n-r)=\frac{2n}{3}\mb\mb\mb(2n)$$
Then Proposition $1$ rewrites into Theorem $4$ as follows.
\begin{Theorem} Assume $6\leq 2n\leq p-1$. Then,
\begin{equation}
A_{2n}=\frac{p^2}{2}\mc\mb(2n)-p\,\mb_{2n}-\frac{p^3}{6}\,\mbt(2n)-\frac{16n^2-12n+5}{24}\,p^3\,\mb_{2n-2}\mpq
\end{equation}
\end{Theorem}
\newtheorem{Remark}{Remark}
\begin{Remark}
The first two terms of Congruence $(10)$ taken modulo $p^3$ are precisely $A_{2n}$ mod $p^3$.
\end{Remark}
We will need a statement similar to Theorem $5$ point II of \cite{LEV3} in order to conclude.
We state it, then prove it.
\begin{Theorem}\hfill\\
(i) Let $n$ be an integer with $4\leq 2n\leq p-7$. Then,
\begin{equation*}
\as_{2n}=-A_{p-1-2n}+p^2\,w_p(1+pw_p)\mb_{p-1-2n}-\frac{p^3\,w_p}{2}\,\mc\mb(p-1-2n)\;\mpq
\end{equation*}
(ii) Case $2n=2$.\\
\begin{equation*}
\as_2=-A_{p-3}+p^2\,w_p(2\mb_{p-3}-\mb_{2p-4})\mpq
\end{equation*}
(iii) Case $2n=p-5$. \\
\begin{equation*}
\as_{p-5}=\frac{1}{8}\Big(\frac{5p^3}{4}-\frac{p}{15}-\frac{p^2}{36}\Big)+\frac{p\,w_p}{120}\Big(\frac{7p^2}{12}-p\,(1+p\,(pB_{p-1})_1)\Big)\mpq
\end{equation*}
(iv) Case $2n=p-3$.\\
\begin{equation*}
\as_{p-3}=\frac{1}{24}\Big(p(2p-1)(5p-2)+p\,w_p\big(2p-11p^2+2p^2(p\,B_{p-1})_1\big)\Big)\mpq
\end{equation*}
$(v)$ Case $2n=p-1$.\\
\begin{equation*}
\as_{p-1}=-1+p\,w_p\Big(3(p\,\mb_{p-1}-1)-p\,\mb_{2(p-1)}-\frac{1}{2}p^2\mb_{p-1}^2\Big)\mpq
\end{equation*}
\end{Theorem}

\textsc{Proof.} Assume first $2\leq 2n\leq p-5$. Recall from Result $9$ that for this range of $2n$, we have:
\begin{equation}
\as_{2n}=p\,\big(2\,\mb_{p-1-2n}-\mb_{2(p-1)-2n}\big)+\frac{p^2}{2}\mt\mc\mb(p+1-2n,p-3)\mpt
\end{equation}
From there, multiply both sides of the latter congruence by $p\,w_p$ and obtain:
\begin{equation}\begin{split}
\as_{2n}=-A_{p-1-2n}+p^2w_p&\big(2\,\mb_{p-1-2n}-\mb_{2(p-1)-2n}\big)\\&+\frac{p^3}{2}\,w_p\,\mt\mc\mb(p+1-2n,p-3)\;\mpq\end{split}
\end{equation}
Suppose first $4\leq 2n\leq p-7$ and apply Result $10$ $(i)$ in order to get point $(i)$ of Theorem $5$. \\
When $2n=2$, the truncated convolution in $(12)$ vanishes and it leads to point $(ii)$ of Theorem $5$.\\
Suppose now $2n=p-5$. We apply point $(iii)$ of Result $10$ to the truncated convolution. It comes:
\begin{equation}\begin{split}
\as_{p-5}=\frac{1}{8}\Big(\frac{5p^3}{4}&-\frac{p}{15}-\frac{p^2}{36}\Big)\\&+p^2w_p\big(2\,\mb_4-\mb_{p+3}\big)+
\frac{pw_p}{2}\Big(\frac{7}{720}p^2+2p\,\mb_{p+3}+2p\,\mb_4(p\,B_{p-1})\Big)\\&\qquad\qquad\qquad\qquad\qquad\qquad\qquad\qquad\qquad
\qquad\qquad\mpq\end{split}
\end{equation}
After simplification and using the fact that $p\,B_{p-1}=-1\mpu$, we obtain point $(iii)$ of Theorem $5$. We verified the congruence on $p=11$ amongst other primes and obtained $A_6^{\star}=2068\;\text{mod}\,11^4$ as holds.\\
For point $(iv)$, we use the expressions for $\as_{p-3}$ modulo $p^3$ (coming from Result $9$ of $\S\,2$) and $A_2$ modulo $p^4$ coming from Proposition $1$ (iii). It yields the result as stated. We verified the congruence on $p=11$ amongst other primes and obtained $A_8^{\star}=5456\,\text{mod}\,11^4$ which is correct.\\
Finally when $2n=p-1$, the multiple harmonic sum is the inverse of the factorial $(p-1)!$. By using the expression for $\as_{p-1}$ modulo $p^3$ as stated in Result $9$ of $\S\,2$, we then derive point $(v)$ of Theorem $5$ which is thus entirely proven.\\

By Theorem $4$ applied with $2n=p-3$ and by using also Theorem $5$ $(ii)$, we have,
\begin{equation}\begin{split}
\frac{p^3}{6}\mbt(p-3)=\as_2&+\frac{p^2}{2}\mc\mb(p-3)\\
&-p\,\mb_{p-3}-\frac{59}{24}\,p^3\,\mb_{p-5}\\
&-p^2\,w_p\,(2\mb_{p-3}-\mb_{2p-4})\mpq
\end{split}\end{equation}

Moreover, from a direct calculation we have: $$\as_2=\frac{1}{2}(H_1^2-H_2),$$ and by Wolstenholme's theorem we know that $p^2|H_1$ (cf Result $4$). Hence,
$$\as_2=-\frac{1}{2}H_2\mpq$$
Then, applying Result $7$ by Sun, we obtain:
\begin{equation}
\as_2=p\big(\mb_{3p-5}-3\,\mb_{2p-4}+3\,\mb_{p-3}\big)+2\,p^3\,\mb_{p-5}\mpq
\end{equation}

Replacing into $(14)$, using Result $11$, Result $12$ $(K2)$ with $k=2$ and $b=p-3$, and simplifying now leads to point $(i)$ of Theorem $1$.

The congruence got tested on $p=11$ and $p=13$, amongst many other primes. In the first case, the second term of the right hand side of the congruence vanishes and in the latter case lots of terms do vanish in the right hand side of the congruence as $13$ is a Wilson prime. Our result is indeed consistent with $\mbt(8)=3\,\text{mod}\,11$ and with $\mbt(10)=2\,\text{mod}\,13$.

Similarly, by applying Theorem $4$ with $2n=p-5$, Theorem $5$ point $(i)$ with $2n=4$, Result $11$ of $\S\,2$ taken from \cite{LEV3} and after simplifying, we get:
\begin{equation}\begin{split}
p^3&\mbt(p-5)=3p^3\big(\mc\mb(p-5)\big)_1-\frac{135}{4}\,p^3\,\mb_{p-7}-6p\mb_{p-5}+6\as_4\\
&+3p^2(1+p\,w_p)\Big\lbrace\Big(2w_p\mb_{p-5}+2(\mb_{2(p-3)}-\mb_{p-5})_1-\mb_{p-3}^2\Big)_0-2w_p\mb_{p-5}\Big\rbrace\\&
\qquad\qquad\qquad\qquad\qquad\qquad\qquad\qquad\qquad\qquad\qquad\qquad\qquad\mpq\end{split}\end{equation}
Moreover, from a direct calculation using the same technique as in \cite{LEV3} and using Bayat's generalization of Wolstenholme's theorem in \cite{BA}, the following congruence holds:
\begin{equation}
\as_4=\frac{1}{8}\big(H_2^2-2\,H_4\big)\mpq
\end{equation}
Further, an application of Sun's result, stated as Result $7$ in $\S\,2$ leads to:
\begin{eqnarray*}
H_2^2&=&4p^2\big(\mb_{3p-5}-3\mb_{2p-4}+3\mb_{p-3}\big)^2\qquad\qquad\;\;\;\;\,\,\mpq\\
H_4&=&-4p\big(\mb_{3p-7}-3\mb_{2p-6}+3\mb_{p-5}\big)-20\,p^3\mb_{p-7}\mpq
\end{eqnarray*}
And so, \begin{equation}\begin{split}
\as_4=&p\big(\mb_{3p-7}-3\mb_{2p-6}+3\mb_{p-5}\big)+\frac{1}{2}\,p^2(\mb_{p-3})_0^2\\&+p^3\mb_{p-3}\big(2\mb_{p-3}-\mb_{2p-4}\big)_1
+5p^3\mb_{p-7}\mpq\end{split}
\end{equation}
Moreover, by Sun's generalization of Kummer's congruences to the modulus $p^2$ in Result $12$ applied with $k=2$ and $b=p-5$, we have
\begin{equation}
\mb_{3p-7}=2\,\mb_{2p-6}-\mb_{p-5}\mpd
\end{equation}
Therefore,
\begin{equation}
\mb_{3p-7}-3\mb_{2p-6}+2\mb_{p-5}=\mb_{p-5}-\mb_{2p-6}\mpd\\
\end{equation}
%In particular,
%\begin{equation}
%(\mb_{3p-7})_1=2(\mb_{2p-6})_1-(\mb_{p-5})_1
%\end{equation}
We thus get point $(ii)$ of Theorem $1$. The congruence got tested on various values of primes using Mathematica.

\section{Wilson's theorem to the modulus $p^4$}

This part is concerned with generalizing Wilson's theorem to the modulus $p^4$. First, we will establish some novel congruence for the generalized harmonic number $H_{p-1,p-1}$ modulo $p^4$. This case is not treated in \cite{SU2}, where a formula for $H_{p-1,k}$ modulo $p^4$ gets only provided for the integers $k$ with $1\leq k\leq p-5$. \\\\
The part is based on Newton's formulas applied with the multiple harmonic sums and the generalized harmonic numbers.\\
In Result $1$ of $\S\,2$, our complex numbers are now chosen to be the first $p-1$ reciprocals of integers. And so $\mcp_k=H_k$ and $\ma_k=\as_k$.
Then, Equality $(1)$ rewrites as:
\begin{equation}
\as_k=\frac{(-1)^{k-1}}{k}\Bigg(H_k+\sum_{r=1}^{k-1}(-1)^r\as_r\,H_{k-r}\Bigg)\qquad\forall k=1,\dots,p-1
\end{equation}
Note, when $k=1$, $\as_1=H_1$, hence the equality still holds. \\\\

In what follows, $k$ is an integer with $k\in\lb p-1,p-3\rb$. \\\\Our first goal will be to compute $H_k$ modulo $p^4$. To that aim, we adapt Sun's proof of \cite{SU2}. Like he does in \cite{SU2} in the modulus $p^3$ case, we use Euler's totient function in order to reduce modulo $p^4$ the generalized harmonic sum to a sum of powers of the first $(p-1)$ integers, which we know how to compute using Bernoulli's formula and reduce further using some novel congruences of Kummer type which Sun develops in \cite{SU1} and \cite{SU2}.

We have:
\begin{equation}
H_k=\sum_{a=1}^k\frac{1}{a^k}=\sum_{a=1}^{p-1}a^{\eu(p^4)-k}\mpq
\end{equation}

Because $\eu(p^4)=p^3(p-1)$ is even and so is $k$, applying Result $3$ yields:
\begin{equation}
H_k=p\,B_{\eu(p^4)-k}+\frac{p^3}{6}(\eu(p^4)-k)(\eu(p^4)-k-1)B_{\euq-k-2}\mpq
\end{equation}

We have,
$$\frac{B_{\eu(p^4)-k}}{\eu(p^4)-k}=\frac{B_{(p-1)(p^3-1)+p-1-k}}{(p-1)(p^3-1)+p-1-k}$$
We must distinguish between $k=p-3$ and $k=p-1$. \\
Suppose first $k=p-3$. We apply $(K3)$:
\begin{equation}\begin{split}
p\,B_{(p-1)(p^3-1)+2}=&p\,\Big((p-1)(p^3-1)+2)\Big)\Big(\binom{p^3-1}{2}\mb_{2(p-1)+2}\\&-(p^3-1)(p^3-3)\mb_{p-1+2}+\binom{p^3-2}{2}(1-p)\mb_2\Big)\mpq
\end{split}\end{equation}
Simplifying $(24)$ leads to:
\begin{equation}
p\,B_{(p-1)(p^3-1)+2}=-\frac{1}{2}\Big(p-3+\frac{11p^3}{2}\Big)\,B_{2p}+(3p-p^2)\Big(-3\mb_{p+1}-\frac{p-1}{4}\Big)\mpq
\end{equation}
Suppose now $k=p-1$. This time, we apply $(S4)$ and get:
\begin{equation}\begin{split}
p\,B_{(p^3-1)(p-1)}=&\frac{13p^3+12p-12}{3}\\&+\frac{19p^3-12}{2}\,pB_{p-1}+(4-7p^3)\,pB_{2(p-1)}+\frac{11p^3-6}{6}\,pB_{3(p-1)}\\
&\qquad\qquad\qquad\qquad\qquad\qquad\qquad\qquad\qquad\qquad\mpq
\end{split}\end{equation}
It remains to deal with the second term of $(23)$. \\
Suppose first $k=p-3$. Then the term reads, where we apply $(S2)$:
\begin{equation}\frac{(p-3)(p-2)}{6}p^2\Big(pB_{(p-1)(p^3-1)}\Big)=p^2\Big(2(p-1)-pB_{p-1}\Big)+\frac{5}{6}p^3\mpq
\end{equation}
And when $k=p-1$, we have $\eu(p^4)-k-2=p^3(p-1)-(p-1)-2$. In particular, $B_{\eu(p^4)-k-2}$ is a $p$-adic integer by Von Staudt's theorem, and so the whole term simply vanishes modulo $p^4$.\\

The results from above get gathered in Theorem $2$ which is thus proven. \\

\noindent All the numbers present in Eq. $(21)$ may now be reduced appropriately. \\
By Result $9$, all of $\as_1,\as_3,\dots,\as_{p-4}$ are divisible by $p^2$. Also, by Result $6$ due to Sun, $p^2$ divides all of $H_1,H_3,\dots,H_{p-4}$ (of course the latter fact was also known since Bayat \cite{BA}). Assume now $k$ is even in $(21)$. Then, modulo $p^4$, only the even indices contribute to the sum of $(21)$, with the addition of $r=1$ and of $r=p-2$ when we deal with $\as_{p-1}$. We are now in a position to provide the multiple harmonic sums $\mathcal{H}_{\lbrace s\rbrace^{2l}=1;p-1}$ with $2\leq 2l\leq p-1$ modulo $p^4$.

\begin{Proposition} Multiple harmonic sums with even indices modulo $p^4$.\\
(i) Assume that $2n\leq p-5$. Then,
\begin{equation}\begin{split}\as_{2n}=-\frac{1}{2n}\Bigg(&-2n\,p\Big(\mb_{3p-3-2n}-3\,\mb_{2p-2-2n}+3\,\mb_{p-1-2n}\Big)\\&
-\binom{2n+2}{3}\,p^3\,\mb_{p-3-2n}\\
&+(2n-2)p^2(2\mb_{p-3}-\mb_{2p-4})(\mb_{2p-2n}-2\mb_{p+1-2n})\\
&+\sum_{r=4}^{2n-2}p^2(2n-r)\bigg(\mb_{2(p-1)-2n+r}-2\mb_{p-1-2n+r}\bigg)\,\mb_{p-1-r}\\
&+\sum_{r=4}^{2n-2}(2n-r)\frac{p^2}{2}\mc\mb(p-1-r)\,p\,\mb_{p-1-2n+r}\\
&-p^3w_p\sum_{r=4}^{2n-2}(2n-r)\mb_{p-1-2n+r}\,\mb_{p-1-r}\Bigg)\;\;\mpq
\end{split}\end{equation}
(ii) Case $2n=p-3$.
\begin{equation}\begin{split}
\as_{p-3}=-\frac{1}{p-3}\Bigg(&\frac{37p^3}{12}-3p^2+\frac{3p}{4}-p^2(pB_{p-1})+3p(p-3)\mb_{p+1}\\
&-\frac{1}{2}\Big(p-3+\frac{11p^3}{2}\Big)B_{2p}\\
&+p^2(p-5)(\mb_{p+3}-2\mb_4)(2\,\mb_{p-3}-\mb_{2p-4})\\
&+2p(\mb_{2p-4}-2\mb_{p-3})\Big(p(2\mb_4-\mb_{p+3})\\&\qquad\qquad\qquad\qquad+\frac{1}{2}\big(\frac{7}{720}p^2+2p\,\mb_{p+3}+2p\mb_4(p\,B_{p-1})\big)\Big)\\
&+\sum_{r=4}^{p-7}p^2(p-3-r)\bigg(\mb_{p+1+r}-2\mb_{r+2}\bigg)\,\mb_{p-1-r}\\
&+\sum_{r=4}^{p-7}(p-3-r)\frac{p^2}{2}\mc\mb(p-1-r)\,p\,\mb_{r+2}\\
&-p^3w_p\sum_{r=4}^{p-7}(p-3-r)\mb_{r+2}\,\mb_{p-1-r}\Bigg)\;\;\mpq
\end{split}\end{equation}
(iii) Case $2n=p-1$.
\begin{equation}\begin{split}
\as_{p-1}=-\frac{1}{p-1}\Bigg(&\frac{13p^3+12p-12}{3}+\frac{19p^3-12}{2}\,pB_{p-1}+(4-7p^3)\,pB_{2(p-1)}\\&+\frac{11p^3-6}{6}\,pB_{3(p-1)}\\
&+p(2\,\mb_{p-3}-\mb_{2p-4})\Big(\big(\frac{1}{2}-3B_{p+1}\big)p-\frac{4}{3}p^2\Big)+\frac{p^3}{2}\mb_{p-3}\\
&+2p(\mb_{2p-4}-2\mb_{p-3})\Big(\frac{p}{12}-\frac{11p^2}{24}+\frac{p^2}{12}(pB_{p-1})_1\Big)\\
&+4p(\mb_{2p-6}-2\mb_{p-5})\Big(p(2\mb_4-\mb_{p+3})\\&\qquad\qquad\qquad\qquad+\frac{1}{2}\big(\frac{7}{720}p^2+2p\,\mb_{p+3}+2p\mb_4(p\,B_{p-1})\big)\Big)\\
&+\sum_{r=4}^{p-7}p^2(p-1-r)\bigg(\mb_{p-1+r}-2\mb_{r}\bigg)\,\mb_{p-1-r}\\
&+\sum_{r=4}^{p-7}(p-1-r)\frac{p^2}{2}\mc\mb(p-1-r)\,p\,\mb_{r}\\
&-p^3w_p\sum_{r=4}^{p-7}(p-1-r)\mb_{r}\,\mb_{p-1-r}\Bigg)\;\;\mpq
\end{split}\end{equation}
\end{Proposition}
\begin{Remark}
Congruence $(28)$ is consistent modulo $p^3$ with Sun's generalization of Kummer congruence $(K2)$ with $k=2$ and $b=p-1-2n$.
Indeed, as none of $\mc\mb(p-1-r)$, $\mb_{p-1-2n+r}$, $\mb_{p-1-r}$, $\mb_{p-3-2n}$ does have its denominator divisible by $p$, the last two terms of $(28)$ vanish modulo $p^3$ and so does the term at the end of the first row. Further, for the range of $r$ that is considered and since $r$ is even, we have by $(K1)$,
$$\mb_{2(p-1)-2n+r}=\mb_{p-1-2n+r}\mpu$$
Thus, we get:
\begin{equation} \begin{split}
\as_{2n}=&p\,\Big(\mb_{3p-3-2n}-3\,\mb_{2p-2-2n}+3\,\mb_{p-1-2n}\Big)\\
&+\frac{p^2}{2n}\sum_{r=2}^{2n-2}(2n-r)\,\mb_{p-1-2n+r}\,\mb_{p-1-r}\;\;\;\mpt\end{split}\end{equation}
Further, modulo $p^3$, the second row of $(31)$ is precisely
$$\frac{p^2}{2}\mt\mc\mb(p+1-2n,p-3)$$
Then, by using Result $9$, we get:
\begin{equation}
p\,\Big(\mb_{3p-3-2n}-3\,\mb_{2p-2-2n}+3\,\mb_{p-1-2n}\Big)=p\Big(2\mb_{p-1-2n}-\mb_{2(p-1)-2n}\Big)\;\;\mpt
\end{equation}
Hence,
\begin{equation}
\mb_{3(p-1)-2n}=2\,\mb_{2(p-1)-2n}-\mb_{p-1-2n}=0\;\mpd
\end{equation}
\end{Remark}
Going back to the modulus $p^4$, since $4\leq p-1-2n+r\leq p-3$, Result $13$ applies and yields:
\begin{equation}
\mb_{2(p-1)-2n+r}=\mb_{p-1-2n+r}-\frac{p}{2}\sum_{a=1}^{p-1}\frac{q_a^2}{a^{2n-r}}\mpd
\end{equation}
Thus, the first and last sums of $(28)$ result in
\begin{equation}\begin{split}\msu+\mst=
&-p^2n(1+pw_p)\mt\mc\mb(p+1-2n,p-3)\\&+p^2(1+pw_p)(2n-2)\mb_{p+1-2n}\mb_{p-3}\\
&-\frac{p^3}{2}\sum_{a=1}^{p-1}\Big(\sum_{i=p+1-2n}^{p-5}\frac{(2n+1)\mb_i+B_i}{a^i}\Big)\frac{q_a^2}{a^{2n}}\;\mpq
\end{split}\end{equation}
We now process the second sum of $(28)$. First and foremost, we do the change of indices $i=p-1-2n+r$.
We get, where we also used Result $11$:
\begin{equation}\begin{split}
\msd=&-\frac{p^3}{2}\Bigg(\sum_{i=p+1-2n}^{p-3}\mc\mb(2(p-1)-2n-i)\mb_i\\&+\sum_{i=p+1-2n}^{p-3}\mc\mb(2(p-1)-2n-i)B_i\Bigg)\\
&-p^3(n-1)\Big(2w_p\mb_{p-3}+2(\mb_{2p-4}-\mb_{p-3})_1\Big)\mb_{p+1-2n}\;\mpq
\end{split}\end{equation}
We first tackle the first sum of $(36)$. If in that sum we consider the full sum instead, starting at $i=2$ and ending at $i=2(p-1)-2n-4$, we get a sum of products of three divided Bernoulli numbers whose indices sum up to $2(p-1)-2n$. Hence, the sum is:
\begin{equation*}\begin{split}
\mbt(2(p-1)-2n)-\sum_{i=2}^{p-1-2n}&\mc\mb(2(p-1)-2n-i)\mb_i\\&-\sum_{i=p-1}^{2(p-1)-2n-4}\mc\mb(2(p-1)-2n-i)\mb_i
\end{split}\end{equation*}
And it is also after change of indices $j=2(p-1)-2n-i$ in the last sum:
$$\mbt(2(p-1)-2n)-\sum_{i=2}^{p-1-2n}\mc\mb(2(p-1)-2n-i)\mb_i-\sum_{j=4}^{p-1-2n}\mc\mb(j)\mb_{2(p-1)-2n-j}$$
Moreover, this sum needs to be worked out modulo $p$ and except when $j=p-1-2n$, we have by $(K1)$,
$$\mb_{2(p-1)-2n-j}=\mb_{p-1-2n-j}\mpu$$
And so the first sum of $(36)$ decomposes modulo $p$ as
\begin{equation}\begin{split}
\sum_{i=p+1-2n}^{p-3}&\mc\mb(2(p-1)-2n-i)\mb_i=\mbt(2(p-1)-2n)-\mc\mb(p-1-2n)\mb_{p-1}\\
&-\sum_{i=2}^{p-1-2n}\mc\mb(2(p-1)-2n-i)\mb_i-\sum_{j=4}^{p-3-2n}\mc\mb(j)\mb_{p-1-2n-j}\;\mpu
\end{split}\end{equation}
As for the second sum of $(36)$, we get:
$$\mb\mb B(2(p-1)-2n)-\sum_{i=2}^{p-1-2n}\mc\mb(2(p-1)-2n-i)B_i-\sum_{j=4}^{p-1-2n}\mc\mb(j)B_{2(p-1)-2n-j}$$
Moreover, except when $j=p-1-2n$, a quick calculation using $(K1)$ shows that
\begin{equation}
B_{2(p-1)-2n-j}=B_{p-1-2n-j}-\mb_{p-1-2n-j}\mpu
\end{equation}
From there, summing both sums of $(36)$ leads to
$$\mbt(2(p-1)-2n)+\mb\mb B(2(p-1)-2n)-\mc\mb(p-1-2n)(\mb_{p-1}+B_{p-1})$$
$$-\sum_{j=4}^{p-1-2n-2}\mc\mb(j)B_{p-1-2n-j}-\sum_{i=2}^{p-1-2n}\mc\mb(2(p-1)-2n-i)(\mb_i+B_i)$$
Further, the first sum of the second row is nothing else than $\mb\mb B(p-1-2n)$. \\
Assume that $2\leq i\leq p-7-2n$. We have, where we work modulo $p^4$:
\begin{equation*}\begin{split}p^3\mc\mb(2(p-1)-2n-i)\equiv p^3\big(2\mc\mb(p-1&-2n-i)+2\mb_{p-1-2n-i}\mb_{p-1}\\&+\mt\mc\mb(p+1-2n-i,p-3)\big)\end{split}\end{equation*}
By using point $(i)$ of Result $10$ and the fact that $p^3(\mb_{p-1}+w_p)=p^2\mpq$, it comes:
$$p^3\mc\mb(2(p-1)-2n-i)\equiv p^3\mc\mb(p-1-2n-i)+2p^2\mb_{p-1-2n-i}-p^3\sum_{a=1}^{p-1}\frac{q_a^2}{a^{2n+i}}$$
By gathering all the contributions and computing further, we obtain Theorem $6$ below. It got tested successfully on several admissible couples $(p,n)$ using a Mathematica program.

\begin{Theorem} Multiple harmonic sums $\as_{2n}$ modulo $p^4$ when $2n\leq p-5$ are congruent to:
\begin{equation}\begin{split}
-\frac{1}{2n}\Bigg(&
\frac{p^3}{2}\Big(\frac{2n-1}{3}(\mbt (2(p-1)-2n)-\mbt (p-1-2n))\\&-\frac{2n}{3}\mbt(p-1-2n)\Big)\\
&-\frac{p^3}{2}\sum_{a=1}^{p-1}\Bigg(\sum_{i=p+1-2n}^{p-5}\frac{(2n+1)\mb_i+B_i}{a^i}+\sum_{i=2}^{p-7-2n}\frac{\mb_i+B_i}{a^i}\Bigg)\frac{q_a^2}{a^{2n}}\\
&-2n\,p\Big(\mb_{3p-3-2n}-3\,\mb_{2p-2-2n}+3\,\mb_{p-1-2n}\Big)-\binom{2n+2}{3}\,p^3\,\mb_{p-3-2n}\\
&+2(n-1)p^2\mb_{p-3}(\mb_{2p-2n}-\mb_{p+1-2n})\\&-p^2n(1+pw_p)\mt\mc\mb(p+1-2n,p-3)\\
%&-p^2n(1+pw_p)\mt\mc\mb(p+1-2n,p-3)%+\mc\mb(p-1-2n)\Big(\frac{p^3}{2}(\mb_{p-1}+B_{p-1})-p^4\mb_{p-1}\Big)\\
%+\mc\mb(p-1-2n)\Big(\frac{p^3}{2}(\mb_{p-1}+B_{p-1}+1)-\frac{2n-1}{2}p^2\Big)\\
&-\frac{2n-1}{2}p^2\mc\mb(p-1-2n)+\frac{p^3}{2}\mc\mb(p-1)(\mb_{p-1-2n}+B_{p-1-2n})\\
&+\Big(\frac{p^3}{2}\big(\mt\mc\mb(4,p-3)+2\mb_2\mb_{p-1}\big)-p^2\mb_2\Big)(\mb_{p-3-2n}+B_{p-3-2n})\\
&+\Big(\frac{p^3}{2}\big(\mt\mc\mb(6,p-3)+2\mb_4\mb_{p-1}+\mb_2^2\big)-p^2\mb_4\Big)(\mb_{p-5-2n}+B_{p-5-2n})\Bigg)
\end{split}\end{equation}
\end{Theorem}

\begin{Remark}
After inspection, the expression of Theorem $6$ is consistent with the modulus $p^3$ case as described in Result $9$. %In order to see this, it is necessary to use Fact $2$ of \cite{LEV1}.
\end{Remark}
Following the same sketch of proof, we can treat the case $2n=p-3$ as well and we obtain the following congruence, where the first three rows have been copied from before.
\begin{Theorem}
\begin{equation}\begin{split}
\as_{p-3}=&-\frac{1}{p-3}\Bigg(\frac{37p^3}{12}-3p^2+\frac{3p}{4}-p^2(pB_{p-1})+3p(p-3)\mb_{p+1}\\&-\frac{1}{2}\Big(p-3+\frac{11p^3}{2}\Big)B_{2p}\\
&+p^2(p-5)(\mb_{p+3}-2\mb_4)(2\,\mb_{p-3}-\mb_{2p-4})\\
&+2p(\mb_{2p-4}-2\mb_{p-3})\Big(p(2\mb_4-\mb_{p+3})\\&+\frac{1}{2}\big(\frac{7}{720}p^2+2p\,\mb_{p+3}+2p\mb_4(p\,B_{p-1})\big)\Big)\\
&-\frac{p^3}{2}\Big(\frac{4}{3}\mbt(p+1)-\mc\mb(p-1)(\mb_2+B_2)-\mc\mb(p-3)(\mb_4+B_4)\\&-\mb_2^2(\mb_{p-3}+B_{p-3})\Big)\\
&+p^2(1+pw_p)\frac{p-3}{2}\Big(2\mb_4\mb_{p-3}-\mt\mc\mb(4,p-3)\Big)\\&-\frac{p^3}{2}\sum_{a=1}^{p-1}\Bigg(\sum_{i=6}^{p-5}\frac{B_i-2\mb_i}{a^i}\Bigg)(aq_a)^2\Bigg)\\
&\qquad\qquad\qquad\qquad\qquad\qquad\qquad\qquad\qquad\qquad\qquad\mpq
\end{split}\end{equation}
\end{Theorem}
Doing the same kind of work when $2n=p-1$ leads to the theorem  below.
\begin{Theorem}$\as_{p-1}$ is congruent modulo $p^4$ to
\begin{equation}\begin{split}
&-\frac{1}{p-1}\Bigg(-\frac{p^3}{3}\mbt (p-1)+\frac{p^3}{4}\Big(w_p\,\mb_{p-3}+(\mb_{2p-4}-\mb_{p-3})_1\Big)\\&-\frac{p^3}{72}\mb_{p-5}+\frac{7p^3}{12}\mb_{p-3}\\
&+\frac{1}{2}p^3\Big(\mc\mb(p-1)\Big)_1+\frac{1-p}{2}(1+pw_p)p^2(2p\mb_{2(p-1)}-p^2\mb_{p-1}^2)_2\\
&-\frac{p^2}{6}(1+pw_p)(\frac{1}{2}\mb_{p-3}+\frac{1}{5}\mb_{p-5})-\frac{p^3}{2}\sum_{a=1}^{p-1}\sum_{i=6}^{p-5}\frac{B_i}{a^i}\,q_a^2\\
&\frac{13p^3+12p-12}{3}+\frac{19p^3-12}{2}\,pB_{p-1}+(4-7p^3)\,pB_{2(p-1)}+\frac{11p^3-6}{6}\,pB_{3(p-1)}\\
&+p(2\,\mb_{p-3}-\mb_{2p-4})\Big(\big(\frac{1}{2}-3B_{p+1}\big)p-\frac{4}{3}p^2\Big)\\
&+2p(\mb_{2p-4}-2\mb_{p-3})\Big(\frac{p}{12}-\frac{11p^2}{24}+\frac{p^2}{12}(pB_{p-1})_1\Big)\\
&+4p(\mb_{2p-6}-2\mb_{p-5})\Big(p(2\mb_4-\mb_{p+3})+\frac{1}{2}\big(\frac{7}{720}p^2+2p\,\mb_{p+3}+2p\mb_4(p\,B_{p-1})\big)\Big)\Bigg)
%&\qquad\qquad\qquad\qquad\qquad\qquad\qquad\qquad\qquad\qquad\qquad\qquad\qquad\qquad\mpq
\end{split}\end{equation}
\end{Theorem}
Before moving further, we will need to study the last sum, say $S$, of the fourth row of $(41)$. In order to do so, we will study a similar sum, namely  $$\mS:=-\frac{p^3}{2}\sum_{a=1}^{p-1}\sum_{i=2}^{p-5}\frac{\mb_{i}}{a^i}q_a^2\mpq$$
We will never compute the sum $\mS$ but rather use it as an artefact for the computation.\\

First and foremost, we will rewrite the sum $\mS$ using Ernvall and Mets\"ankyl\"a's result. According to their result, we have for those indicies $i$ with $2\leq i\leq p-5$:
\begin{equation}
(\mb_{p-1+p-1-i}-\mb_{p-1-i})_1=-\frac{1}{2}\sum_{a=1}^{p-1}\frac{q_a^2}{a^i}\mpu
\end{equation}
Then, the sum of interest rewrites as
\begin{equation}
\mS=p^3\sum_{i=2}^{p-5}(\mb_{p-1+p-1-i}-\mb_{p-1-i})_1(\mb_i)_0\mpq
\end{equation}

It comes:
\begin{equation}\begin{split}
\mS=&p^3\sum_{i=2}^{p-3}(\mb_{2(p-1)-i})_1(\mb_i)_0-p^3\sum_{i=2}^{p-3}(\mb_{p-1-i})_1(\mb_i)_0\\
&+p^3\mb_{p-3}\big((\mb_2)_1-(\mb_{p+1})_1\big)\qquad\qquad\mpq
\end{split}\end{equation}

Next, it will be useful to calculate $(\mc\mb(p-1))_1$ and $(\mc\mb(2(p-1))_1$. We have:
$$(\mc\mb(p-1))_1=2\sum_{j=2}^{p-3}(\mb_j)_1(\mb_{p-1-j})_0+\bigg(\sum_{j=2}^{p-3}(\mb_j)_0(\mb_{p-1-j})_0\bigg)_1\;\;\mpu$$
Write the Hensel expansion of $\mb_{p-1}$ as: $$\mb_{p-1}=\frac{1}{p}+(\mb_{p-1})_0+p(\mb_{p-1})_1+p^2(\mb_{p-1})_2+...$$
Then we also have:
\begin{equation*}\begin{split}(\mc\mb(2(p-1)))_1&=2\sum_{j=2}^{p-3}(\mb_j)_1(\mb_{2(p-1)-j})_0+2\sum_{j=2}^{p-3}(\mb_j)_0(\mb_{2(p-1)-j})_1\\
&+2(\mb_{p-1})_0(\mb_{p-1})_1+2(\mb_{p-1})_2\\&+\Big\lbrace 2\sum_{j=2}^{p-3}(\mb_j)_0(\mb_{p-1-j})_0+(\mb_{p-1})_0^2+2(\mb_{p-1})_1+\big(2(\mb_{p-1})_0\big)_1\Big\rbrace_1\\&\qquad\qquad\qquad\qquad\qquad\qquad
\qquad\qquad\qquad\qquad\qquad\qquad\;\mpu
\end{split}\end{equation*}
The latter two congruences now allow to write:
\begin{equation}\begin{split}
\mS=&\frac{p^3}{2}\big((\mc\mb(2(p-1)))_1-2(\mc\mb(p-1))_1\big)-p^3(\mb_{p-1})_1(\mb_{p-1})_0-p^3(\mb_{p-1})_2\\
&-\frac{p^3}{2}\Big(\big[2(\mc\mb(p-1))_0\big]_0+\big[(\mb_{p-1})_0^2\big]_0+\big[2(\mb_{p-1})_1\big]_0+\big(2(\mb_{p-1})_0\big)_1\Big)_1\\
&-\frac{p^3}{2}\Big(\big[2(\mc\mb(p-1))_0\big]_1+\big[(\mb_{p-1})_0^2\big]_1+\big[2(\mb_{p-1})_1\big]_1\Big)\\
&+p^3\mb_{p-3}\big((\mb_2)_1-(\mb_{p+1})_1\big)\qquad\qquad\qquad\qquad\qquad\qquad\;\mpq
\end{split}\end{equation}
%\mS=\frac{p^3}{2}\big(\mc\mb(2(p-1))-\mc\mb(p-1)\big)_1-\frac{p^3}{2}\big(\mc\mb(p-1)\big)_1-\frac{p^3}{2}(\mb_{p-1}^2)_1
We move on to studying the differential term of $(45)$. Say,
$$\md:=\fcmcd-2\fcmcu$$
The key here will be to apply Miki's identity from \cite{MI} in order to make the term of interest of $(41)$ appear. The process will be similar to the one developed in \cite{LEV3} while we were studying some truncated convolutions of divided Bernoulli numbers. \\\\
By applying Miki's identity, we have:
\begin{eqnarray*}
\mc\mb(2(p-1))&=&\sum_{i=2}^{2(p-1)-2}\binom{2(p-1)}{i}\mb_i\mb_{2(p-1)-i}+2\mh_{2(p-1)}\mb_{2(p-1)}\\
&=&2\sum_{i=2}^{p-3}\binom{2(p-1)}{i}\mb_i\mb_{2(p-1)-i}+\binom{2(p-1)}{p-1}\mb_{p-1}^2\\&&\qquad\qquad\qquad\qquad\qquad+2\mh_{2(p-1)}\mb_{2(p-1)}
\end{eqnarray*}
Further, by the Chu-Vandermonde convolution,
$$\binom{p-1+p-1}{i}=\sum_{k=0}^i\binom{p-1}{k}\binom{p-1}{i-k}\qquad\qquad(CVC)$$
Then, \begin{eqnarray*}\binom{2(p-1)}{i}&=&\sum_{k=0}^i(-1)^k(-1)^{i-k}\mpu\\
&=&i+1\qquad\qquad\;\;\;\;\;\mpu\end{eqnarray*}
Then, we get, where we used Ernvall and Mets\"ankyl\"a's congruence again:
%\begin{eqnarray*}
%\fcmcd&=&
%\Bigg\lbrace\big(\mc\mb(p-1)\big)_0+\bigg[\binom{2(p-1)}{p-1}\mb_{p-1}^2-2\mh_{2(p-1)}\mb_{2(p-1)}\bigg]_0\Bigg\rbrace_1\\
%&&+\Bigg\lbrace %2\bigg(\frac{1}{2}\big(\mc\mb(p-1)\big)_0+\big(\mb_{p-3}B_2\big)_0\bigg)_0\Bigg\rbrace_1+\Big(2\binom{2(p-1)}{p-3}\mb_{p-3}\mb_{p+1}\Big)_1\\
%\\&&+2\sum_{i=2}^{p-5}\binom{2(p-1)}{i}_1\mb_i\mb_{p-1-i}+2\sum_{i=2}^{p-5}(i+1)(\mb_i)_1\mb_{p-1-i}
%\\&&+
%2\sum_{i=2}^{p-5}(i+1)\mb_i(\mb_{p-1-i})_1-\sum_{i=2}^{p-5}(i+1)\mb_i\sum_{a=1}^{p-1}\frac{q_a^2}{a^i}\\
%&&+\big(\mc\mb(p-1)\big)_0+2\Bigg(\frac{\mc\mb(p-1)}{2}\Bigg)_1+2\big(2\mb_2\mb_{p-3})_1-2\big(\mb_{p-3}\mb_2\big)_0\\
%&&+2\bigg(\frac{1}{2}\big(\mc\mb(p-1)\big)_0+\big(\mb_{p-3}B_2\big)_0\bigg)_1\\
%&&+\Big[\binom{2(p-1)}{p-1}\mb_{p-1}^2-2\mh_{2(p-1)}\mb_{2(p-1)}\Big]_1
%\qquad\qquad\mpu
%\end{eqnarray*}
\begin{equation*}\begin{split}
&\fcmcd=\Big[2\Big(\frac{1}{2}\mc\mb(p-1)\Big)_0\Big]_1+
2\Big(\sum_{i=2}^{p-3}(i+1)(\mb_i)_0(\mb_{p-1-i})_0\Big)_1\\
&+2\sum_{i=2}^{p-3}\binom{2(p-1)}{i}_1\mb_i\mb_{p-1-i}+2\sum_{i=2}^{p-3}(i+1)(\mb_i)_1\mb_{p-1-i}
\\&+
2\sum_{i=2}^{p-5}(i+1)\mb_i(\mb_{p-1-i})_1-\sum_{i=2}^{p-5}(i+1)\mb_i\sum_{a=1}^{p-1}\frac{q_a^2}{a^i}-4(\mb_{p-3})_0(\mb_{p+1})_1\\
&\Big[\binom{2(p-1)}{p-1}\mb_{p-1}^2\Big]_1+\Big[2\mh_{2(p-1)}\mb_{2(p-1)}\Big]_1+\Big[p-1+\big(2+2(\mb_{p-1})_0\big)_0\Big]_2\\
&+\Bigg\lbrace\big(\mc\mb(p-1)\big)_0+\Big[p-1+\big(2+2(\mb_{p-1})_0\big)_0\Big]_1\\&\qquad+\Big[\binom{2(p-1)}{p-1}\mb_{p-1}^2\Big]_0+\Big[2\mh_{2(p-1)}\mb_{2(p-1)}\Big]_0\Bigg\rbrace_1
\qquad\mpu
\end{split}\end{equation*}
\noindent Furthermore, by Wolstenholme and Von-Staudt's theorems, we know that $$v_p(\mh_{p-1}\mb_{p-1})\geq 1$$
Then, we simply have,
\begin{multline*}
\fcmcu=\sum_{i=2}^{p-3}\binom{p-1}{i}_1\mb_i\mb_{p-1-i}+\sum_{i=2}^{p-3}(\mb_i)_1\mb_{p-1-i}+\sum_{i=2}^{p-3}\mb_i(\mb_{p-1-i})_1\\
+\big[2\mh_{p-1}\mb_{p-1}\big]_1+\Big(\sum_{i=2}^{p-3}(\mb_i)_0(\mb_{p-1-i})_0\Big)_1\mpu
\end{multline*}
We are now in a position to compute $\md$.
We obtain:
%\begin{eqnarray}
%\md&=&2\sum_{i=2}^{p-5}\mb_i\mb_{p-1-i}\Bigg(\binom{2(p-1)}{i}_1-\binom{p-1}{i}_1\Bigg)\\
%&&+2\sum_{i=2}^{p-5}i\,(\mb_i)_1\mb_{p-1-i}+2\sum_{i=2}^{p-5}i\,\mb_i(\mb_{p-1-i})_1\\
%&&-\sum_{i=2}^{p-5}\sum_{a=1}^{p-1}\frac{\mb_i}{a^i}q_a^2-\sum_{i=2}^{p-5}\sum_{a=1}^{p-1}\frac{B_i}{a^i}q_a^2\\
%&&+2\Big(\binom{2(p-1)}{p-3}_1-\binom{p-1}{p-3}_1\Big)\mb_2\mb_{p-3}-\big(2(2\mb_{p-3}\mb_2)_0\big)_1\\
%&&+2(\mb_{p-3})_0\big((\mb_2)_1-2(\mb_{p+1})_1\big)-2(\mb_2)_0\big((\mb_{p-3})_0+(\mb_{p-3})_1\big)\\
%%&&-6\mb_2(\mb_{p-3})_1-2\mb_{p-3}\big(2(\mb_{p+1})_1+(\mb_2)_1\big)\\
%&&+\Big[\binom{2(p-1)}{p-1}\mb_{p-1}^2-2\mh_{2(p-1)}\mb_{2(p-1)}\Big]_1-2\big[-2(\mh_{p-1}\mb_{p-1}\big]_1\\
%&&+\Bigg\lbrace\big(\mc\mb(p-1)\big)_0+\bigg[\binom{2(p-1)}{p-1}\mb_{p-1}^2-2\mh_{2(p-1)}\mb_{2(p-1)}\bigg]_0\Bigg\rbrace_1\\
%&&+\Bigg\lbrace %2\bigg(\frac{1}{2}\big(\mc\mb(p-1)\big)_0+\big(\mb_{p-3}B_2\big)_0\bigg)_0\Bigg\rbrace_1\\&&+2\bigg(\frac{1}{2}\big(\mc\mb(p-1)\big)_0
%+\big(\mb_{p-3}B_2\big)_0\bigg)_1+\big(\mc\mb(p-1)\big)_0+2\Bigg(\frac{\mc\mb(p-1)}{2}\Bigg)_1\notag\\&&\qquad\qquad\qquad\qquad\qquad\qquad\qquad
%\qquad\qquad\qquad\qquad\qquad\mpu
%\end{eqnarray}
\begin{eqnarray}
\md&=&2\sum_{i=2}^{p-3}\mb_i\mb_{p-1-i}\Bigg(\binom{2(p-1)}{i}_1-\binom{p-1}{i}_1\Bigg)\\
&&\nts\nts\nts\nts\nts\nts\nts\nts\nts+2\sum_{i=2}^{p-3}i\,(\mb_i)_1\mb_{p-1-i}+2\sum_{i=2}^{p-3}i\,\mb_i(\mb_{p-1-i})_1+2\Big(\sum_{i=2}^{p-3}i(\mb_i)_0(\mb_{p-1-i})_0\Big)_1\\
&&\nts\nts\nts\nts\nts\nts\nts\nts\nts+2\Bigg[\Big(-\frac{1}{2}\mc\mb(p-1)\Big)_0+\Big(\mc\mb(p-1)\Big)_0\Bigg]_1+\Bigg[2\Big(\frac{1}{2}\mc\mb(p-1)\Big)_0\Bigg]_1\\
&&\nts\nts\nts\nts\nts\nts\nts\nts\nts+\,4\,\mb_{p-3}\big((\mb_2)_1-(\mb_{p+1})_1\big)\\
&&\nts\nts\nts\nts\nts\nts\nts\nts\nts-\sum_{i=2}^{p-5}\sum_{a=1}^{p-1}\frac{\mb_i}{a^i}q_a^2-\sum_{i=2}^{p-5}\sum_{a=1}^{p-1}\frac{B_i}{a^i}q_a^2\\
&&\nts\nts\nts\nts\nts\nts\nts\nts\nts+\big[2\mh_{2(p-1)}\mb_{2(p-1)}\big]_1-2\big[2\mh_{p-1}\mb_{p-1}\big]_1\\
&&\nts\nts\nts\nts\nts\nts\nts\nts\nts+\Big[\binom{2(p-1)}{p-1}\mb_{p-1}^2\Big]_1+\Big[p-1+\big(2+2(\mb_{p-1})_0\big)_0\Big]_2\\
&&\nts\nts\nts\nts\nts\nts\nts\nts\nts+\Bigg\lbrace\big(\mc\mb(p-1)\big)_0+\Big[p-1+\big(2+2(\mb_{p-1})_0\big)_0\Big]_1\\&&\qquad
\Big[\binom{2(p-1)}{p-1}\mb_{p-1}^2\Big]_0+\Big[2\mh_{2(p-1)}\mb_{2(p-1)}\Big]_0\Bigg\rbrace_1\qquad\mpu\notag
\end{eqnarray}
Recall from $(45)$ that
\begin{equation}\begin{split}
\mS =\frac{p^3}{2}\md+p^3\mb_{p-3}\big((\mb_2)_1-(\mb_{p+1})_1\big)-p^3(\mb_{p-1})_1(\mb_{p-1})_0-p^3(\mb_{p-1})_2&-\frac{p^3}{2}\mo_0\\
&\mpq
\end{split}\end{equation}
where $\mo_0$ denotes the contribution from the terms whose writing has been omitted. Namely,
$$\mo_0=\bigg\lbrace 2\big(2p\mb_{2(p-1)}-p^2\mb_{p-1}^2\big)_2+(\mb_{p-1})_0^2+2(\mb_{p-1})_1+\big(2(\mb_{p-1})_0\big)_1\bigg\rbrace_1$$
It follows that:
\begin{equation}\begin{split}
-\frac{p^3}{2}\sum_{i=2}^{p-5}\sum_{a=1}^{p-1}\frac{B_i}{a^i}q_a^2&=p^3(\mb_{p-1})_1(\mb_{p-1})_0\!+\!p^3(\mb_{p-1})_2-\!3p^3\mb_{p-3}\big((\mb_2)_1-(\mb_{p+1})_1\big)\\
&+\frac{p^3}{2}\mo_0-\frac{p^3}{2}\Big((46)+(47)+(48)+(51)+(52)+(53)\Big)\\
&\qquad\qquad\qquad\qquad\qquad\qquad\qquad\qquad\qquad\qquad\qquad\mpq
\end{split}\end{equation}

We first inspect $(47)$ and $(48)$. The sum in $(47)$ is precisely
$$2\Big(\frac{p-1}{2}\mc\mb(p-1)\Big)_1$$
Moreover, we have:
\begin{eqnarray}
\nts\nts2\Big(\frac{p-1}{2}\mc\mb(p-1)\Big)_1\nts\nts&=&\nts\nts-\big(\mc\mb(p-1)\big)_1+\Big(1+2\left(\begin{array}{l}-\frac{1}{2}\end{array}\right)_1\Big)\big(\mc\mb(p-1)\big)_0\\
&&+2\Big(\left(\begin{array}{l}-\frac{1}{2}\end{array}\right)_0\big(\mc\mb(p-1)\big)_0\Big)_1\;\mpu
\end{eqnarray}
Adding $(48)$ to $(57)$ modulo $p$ leads to:
$$2\Bigg[\Big(-\frac{1}{2}\Big)_0\Big(\mc\mb(p-1)\Big)_0+\Big(-\frac{1}{2}\mc\mb(p-1)\Big)_0+\Big(\mc\mb(p-1)\Big)_0\Bigg]_1-2,$$
which is nothing else than:
$$2\Bigg[\Big(-\frac{1}{2}\Big)_0\Big(\mc\mb(p-1)\Big)_0\Bigg]_1+2\Bigg[2\Big(-\frac{1}{2}\mc\mb(p-1)\Big)_0\Bigg]_1$$
Also, by forthcoming Corollary $2$, the second term of the right hand side of $(56)$ cancels modulo $p$.
Whence,
\begin{equation}\begin{split}
-\frac{p^3}{2}\big((47)+(48)\big)&=\frac{p^3}{2}\big(\mc\mb(p-1)\big)_1-\frac{p^3}{2}\Bigg[2\Big(\frac{1}{2}\mc\mb(p-1)\Big)_0\Bigg]_1\\
&-p^3\Bigg[\Big(-\frac{1}{2}\Big)_0\Big(\mc\mb(p-1)\Big)_0\Bigg]_1-p^3\Bigg[2\Big(-\frac{1}{2}\mc\mb(p-1)\Big)_0\Bigg]_1\\&
\qquad\qquad\qquad\qquad\qquad\qquad\qquad\qquad\qquad\;\mpq
\end{split}\end{equation}
and $\big(\mc\mb(p-1)\big)_0=\big(2p\mb_{2(p-1)}-p^2\mb_{p-1}^2\big)_2$, cf. Result $11$. \\

%$$2\Bigg[-\Big(\frac{1}{2}\mc\mb(p-1)\Big)_0\Bigg]_1,$$
%which is obviously simply $-2$. And so $-\frac{p^3}{2}\big((47)+(48)\big)$ rewrites as:
%\begin{equation}
%-\frac{p^3}{2}\big((47)+(48)\big)=p^3+\frac{p^3}{2}\big(\mc\mb(p-1)\big)_1
%-\frac{p^3}{2}\Big(1+2\left(\begin{array}{l}-\frac{1}{2}\end{array}\right)_1\Big)\big(\mc\mb(p-1)\big)_0
%\end{equation}
We now deal with $(46)$. To that aim, we go back to $(CVC)$ and investigate the equality one $p$ power further, that is modulo $p^2$. We will use the following lemma.
\newtheorem{Lemma}{Lemma}
\begin{Lemma}
$\binom{p-1}{k}=(-1)^k(1-p\mh_k)\mpd$
\end{Lemma}
\textsc{Proof.} This is Proposition $6$ of \cite{LEV3}. \\\\
In what follows, $i$ denotes an even integer chosen such that $2\leq i\leq p-3$
We thus have (with the convention that $\mh_0=0$):
\begin{equation}
\binom{2(p-1)}{i}=\sum_{k=0}^i(1-p\,\mh_k)(1-p\,\mh_{i-k})\mpd
\end{equation}
It follows that:
\begin{equation}
\binom{2(p-1)}{i}_1=-2\,\sum_{k=1}^i(\mh_k)_0\;\;\mpu
\end{equation}
Then, we have:
\begin{eqnarray}
\binom{2(p-1)}{i}_1-\binom{p-1}{i}_1&=&(\mh_i)_0-2\sum_{k=0}^{i}(\mh_k)_0\;\mpu
\end{eqnarray}
But, by the lemma, $$p\mh_k=1-(-1)^k\binom{p-1}{k}\;\mpd$$
Then, we have:
\begin{eqnarray}
p\sum_{k=0}^i\mh_k&=&i+1-\sum_{k=0}^i(-1)^k\binom{p-1}{k}\qquad\mpd\\
&=&i+1-\binom{p-2}{i}\qquad\qquad\qquad\mpd
\end{eqnarray}
Then, $(61)$ rewrites \textbf{modulo $p$} as:
\begin{eqnarray}
\binom{2(p-1)}{i}_1-\binom{p-1}{i}_1&\equiv &-\frac{2i+1+\binom{p-1}{i}-2\binom{p-2}{i}}{p}\\
&\equiv &-\frac{2i+1+\binom{p-2}{i-1}-\binom{p-2}{i}}{p}\\
&\equiv &-\frac{2i+1-\frac{p-1-2i}{p-1}\binom{p-1}{i}}{p}\\
&\equiv &\Big[\frac{p-1-2i}{p-1}\binom{p-1}{i}-(2i+1)\Big]_1\\
&\equiv &(2i+1)\binom{p-1}{i}_1+2i
\end{eqnarray}
where $(65)$ and $(66)$ are derived using some classical identities on binomial coefficients.
We thus get:
\begin{equation}
-\frac{p^3}{2}(46)=p^3\big(\mc\mb(p-1)\big)_0\;\mpq
\end{equation}
%&=&-(\mh_i)_0+\bigg[\frac{2}{p}\Bigg(\sum_{k=1}^{i-1}(-1)^k\binom{p-1}{k}-(i-1)\!\!\Bigg)\bigg]_0\notag\\
%&&\\
%&=&-(\mh_i)_0+\Bigg[\frac{2\Big((-1)^{i-1}\binom{p-2}{i-1}-i\Big)}{p}\Bigg]_0\\
%&=&\bigg[\frac{\bn{p-1}{i}-1-2\bn{p-2}{i-1}-2i}{p}\bigg]_0
%\end{eqnarray}
%\begin{equation}
%-\frac{p^3}{2}\big[(46)+(49)\big]=p^3\mc\mb(p-1)
%\end{equation}
%In particular, we have:
%\begin{equation}
%\bn{2(p-1}{p-3}_1-\bn{p-1}{p-3}_1=4
%\end{equation}
%as Wolstenholme's theorem forces $(\mh_{p-3})_0=1$.
We will now investigate the binomial term of $(52)$. First and foremost, from
$$\binom{2(p-1)}{p-1}=\sum_{j=0}^{p-1}\binom{p-1}{j}^2,$$
we see that:
\begin{equation}\binom{2(p-1)}{p-1}=0\mpu\end{equation}
Moreover, the binomial version of Wolstenholme's theorem of \cite{WO} asserts that:
$$\bn{2p-1}{p-1}=1\mpt$$
Then,
\begin{eqnarray}
\bn{2p-2}{p-1}&=&\frac{p\;\bn{2p-1}{p-1}}{2p-1}\\
&=&-4p^3-2p^2-p\mpq
\end{eqnarray}
It follows that,
\begin{equation*}\begin{split}
\Bigg[\bn{2(p-1)}{p-1}\mb_{p-1}^2\Bigg]_1=\Bigg[-\Big(\frac{1}{p}+(\mb_{p-1})_0+p(\mb_{p-1}&)_1+p^2(\mb_{p-1})_2\Big)^2\\&
\times (4p^3+2p^2+p)\Bigg]_1\end{split}\end{equation*}
We thus have:
\begin{equation}\begin{split}
-\frac{p^3}{2}(52)=\frac{p^3}{2}\bigg(&\Big(3+2(\mb_{p-1})_0\Big)_1+5+4(\mb_{p-1})_0+(\mb_{p-1})_0^2+2(\mb_{p-1})_1\\
&-\Big[p-1+\big(2+2(\mb_{p-1})_0\big)_0\Big]_2
\bigg)\mpq\\
&
%\mb_{p-1})_1+\frac{p^3}{2}+\frac{p^3}{2}\Big((\mb_{p-1})_0^2+\big[2(\mb_{p-1})_0\big]_1
\end{split}\end{equation}
Also, we have:
\begin{equation}
\Bigg[\bn{2(p-1)}{p-1}\mb_{p-1}^2\Bigg]_0=-3-2(\mb_{p-1})_0\;\mpu
\end{equation}

%We will show the following lemma.
%\begin{Lemma}
%\begin{equation}\Bigg[\bn{2(p-1)}{p-1}\mb_{p-1}^2\Bigg]_1=-4(p^2\mb_{p-1}^2)_0-2(p^2\mb_{p-1}^2)_1-(p^2\mb_{p-1}^2)_2\end{equation}
%\end{Lemma}
%\textsc{Proof of Lemma}.
%By Von Staudt-Clausen's theorem, we know that $p^2$ divides the denominator of $\mb_{p-1}^2$, hence we need to know $\bn{2(p-1)}{p-1}$ modulo %$p^4$. The binomial version of Wolstenholme's theorem of \cite{WO} asserts that
%$$\bn{2p-1}{p-1}=1\mpt$$
%Then,
%\begin{eqnarray}
%\bn{2p-2}{p-1}&=&\frac{p\;\bn{2p-1}{p-1}}{2p-1}\\
%&=&-4p^3-2p^2-p\mpq
%\end{eqnarray}
%From there, the lemma follows.
It remains to tackle $(51)$. To that aim, we write using Wolstenholme's theorem:
$$\mh_{2(p-1)}=p^2u+\frac{1}{p}+1+pv+p^2w\;\mpt,$$
for some $p$-residues $u$, $v$ and $w$. \\
Moreover, the sum $$\Sigma_p:=\frac{1}{p+1}+\dots +\frac{1}{p+(p-2)}$$ verifies the following congruences:
\begin{eqnarray}
\Sigma_p&=&H_1-\frac{1}{p-1}-p\bigg(H_2-\frac{1}{(p-1)^2}\bigg)+p^2\bigg(H_3-\frac{1}{(p-1)^3}\bigg)\,\nts\nts\mpt\\
&=&p^2(4-B_{p-3})+2p+1\qquad\qquad\qquad\qquad\qquad\qquad\;\;\;\;\;\;\,\!\mpt
\end{eqnarray}
where we used Results $5$ and $6(a)$ of Sun (cf $\S\,2$) in order to derive Congruence $(76)$. \\
Then we have, $$v=2\qquad\text{and}\qquad w=4+3(\mb_{p-3})_0\mpu$$ And since $H_1=p^2\mb_{p-3}\mpt$, we have $$u=(\mb_{p-3})_0$$
From there,
$$\begin{array}{l}\mh_{2(p-1)}\mb_{2(p-1)}=\Big(\frac{1}{p}+1+2p+4(1+(\mb_{p-3})_0)p^2\Big)\Big(\big(\frac{1}{2}\big)_0\frac{1}{p}+
\big(\mb_{2(p-1)}\big)_0\\\qquad\qquad\qquad\qquad\qquad\qquad\qquad\qquad\;\;\;+p\big(\mb_{2(p-1)})_1+p^2\big(\mb_{2(p-1)})_2\Big)\;\mpd\end{array}$$
\begin{Lemma} Let $p$ be an odd prime. Then, the $p$-adic expansions of $\frac{1}{2}$ and $-\frac{1}{2}$ are respectively:
\begin{eqnarray*}\frac{1}{2}&=&\frac{p+1}{2}+\sum_{k=1}^{\infty}\frac{p-1}{2}p^k\\&&\\
-\frac{1}{2}&=&\sum_{k=0}^{\infty}\frac{p-1}{2}p^k\end{eqnarray*}
\end{Lemma}
\newtheorem{Corollary}{Corollary}
\begin{Corollary}
$$2\Big(\frac{1}{2}\Big)_0=p+1$$
\end{Corollary}
\begin{Corollary}
$$2\Big(-\frac{1}{2}\Big)_1=p-1$$
\end{Corollary}
\textsc{Proof of Lemma.}\\

\noindent Consider the polynomial $f(X)=2X-1\in\mathbb{Z}_p[X]$. The $p$-residue $\frac{p+1}{2}$ satisfies to
$$\begin{array}{cc}
f\Big(\frac{p+1}{2}\Big)&\in\,p\mathbb{Z}_p\\
f'\Big(\frac{p+1}{2}\Big)&\not\in\,p\mathbb{Z}_p
\end{array}$$
Then, by Hensel's lemma this residue lifts to a unique root $$y=\frac{p+1}{2}+\sum_{i=1}^{\infty}y_ip^i\in\mathbb{Z}_p$$ The expansion of the root is then built inductively by applying Hensel's algorithm. First, $f(\frac{p+1}{2}+y_1p)\in\,p^2\mathbb{Z}_p$ implies immediately $y_1=\frac{p-1}{2}$. \\Next,$f(\frac{p+1}{2}+\sum_{i=1}^l\frac{p-1}{2}p^i+y_{l+1}p^{l+1})\in\,p^{l+2}\mathbb{Z}_p$ implies $y_{l+1}=\frac{p-1}{2}$.\\
From there, the $p$-adic expansion of $-\frac{1}{2}$ can be deduced easily. \hfill $\square$\\

By applying Corollary $1$, we thus get:
\begin{equation}\begin{split}
2\mh_{2(p-1)}\mb_{2(p-1)}=&\frac{1}{p^2}+2\Big(1+\big(\mb_{2(p-1)}\big)_0\Big)\frac{1}{p}+1\\&+2\Big(1+\big(\mb_{2(p-1)}\big)_0+\big(\mb_{2(p-1)}\big)_1\Big)\\
&+p\bigg(2\Big(1+\big(\mb_{2(p-1)}\big)_1+\big(\mb_{2(p-1)}\big)_2\Big)\\&+4\Big(1+\big(\mb_{p-3}\big)_0+\big(\mb_{2(p-1)}\big)_0\Big)\bigg)\mpd
\end{split}\end{equation}

We thus have,
\begin{equation}\begin{split}
\Big[2\mh_{2(p-1)}\mb_{2(p-1)}\Big]_0=\!\Big[2\Big(1+\big(\mb_{2(p-1)}\big)_0\Big)\Big]_1\!+1+2\Big(1+(\mb_{2(p-1)})_0+&(\mb_{2(p-1)})_1\Big)\\&\;\;\mpu
\end{split}\end{equation}
and
\begin{equation}\begin{split}
\Big[2\mh_{2(p-1)}\mb_{2(p-1)}\Big]_1=&\Bigg\lbrace\Big[2\Big(1+\big(\mb_{2(p-1)}\big)_0\Big)\Big]_1
+1\\&+2\Big(1+\big(\mb_{2(p-1)}\big)_0+\big(\mb_{2(p-1)}\big)_1\Big)\Bigg\rbrace_1\\
&+\!\Big[2\Big(1+\big(\mb_{2(p-1)}\big)_0\Big)\Big]_2+\!2\Big(1+\big(\mb_{2(p-1)}\big)_1+\big(\mb_{2(p-1)}\big)_2\Big)\\
&+4\Big(1+\big(\mb_{p-3}\big)_0+\big(\mb_{2(p-1)}\big)_0\Big)\\
&\qquad\qquad\qquad\qquad\qquad\qquad\qquad\qquad\qquad\nts\nts\nts\;\mpu
\end{split}\end{equation}
We also have:
\begin{equation}
-2\Big[2\mh_{p-1}\mb_{p-1}\Big]_1=-4\big(\mb_{p-3}\big)_0\;\qquad\qquad\qquad\qquad\nts\mpu
\end{equation}
It follows that,
\begin{equation}\begin{split}
-\frac{p^3}{2}(51)=&\!-\frac{p^3}{2}\Bigg\lbrace\!\Big[2\Big(1\!+\!\big(\mb_{2(p-1)}\big)_0\Big)\Big]_1
+1+2\Big(1\!+\!\big(\mb_{2(p-1)}\big)_0+\big(\mb_{2(p-1)}\big)_1\Big)\Bigg\rbrace_1\\
&-\frac{p^3}{2}\Big[2\Big(1+\big(\mb_{2(p-1)}\big)_0\Big)\Big]_2-p^3\Big(1+\big(\mb_{2(p-1)}\big)_1+\big(\mb_{2(p-1)}\big)_2\\
&\qquad\qquad\qquad\qquad\qquad\qquad\qquad+2\big(1+\big(\mb_{2(p-1)}\big)_0\big)\Big)\;\mpq
\end{split}\end{equation}
We are ready to write the final form of $(55)$ under the following theorem.
\begin{Theorem}
\begin{equation}
-\frac{p^3}{2}\sum_{i=2}^{p-5}\sum_{a=1}^{p-1}\frac{B_i}{a^i}q_a^2=\frac{p^3}{2}\big(\mc\mb(p-1)\big)_1+\frac{p^3}{2}\mo\;\;\mpq
\end{equation}
where the other term $\mo$ gets provided in Appendix A.
\end{Theorem}
We will now see that Theorem $3$ can be deduced from the conjunction of Theorems $4$, $8$ and $9$. \\
Starting again from Theorem $8$, we may write:
\begin{equation}
\as_{p-1}=-\frac{p^3}{3}\mbt(p-1)+p^3(\mc\mb(p-1))_1+\frac{p^3}{2}\mo+\mo^{'}\mpq
\end{equation}
with $\mo^{'}$ the adequate other terms arising from the right hand side of the congruence of Theorem $8$ and $\mo$ the residue sum, as given in Appendix A.
Now apply Theorem $4$ with $2n=p-1$ and obtain:
\begin{equation}\begin{split}
-\frac{p^3}{3}\mbt(p-1)\!+\!p^3(\mc\mb(p-1))_1\!=\!2A_{p-1}\!-p^2(\mc\mb(p-1))_0\!+\!2p\mb_{p-1}&\!+\!\frac{5}{4}p^3\mb_{p-3}\\
&\mpq
\end{split}\end{equation}
And so,
\begin{equation}
\as_{p-1}-2A_{p-1}=-p^2(\mc\mb(p-1))_0+2p\,\mb_{p-1}+\frac{5}{4}\,p^3\,\mb_{p-3}+\frac{p^3}{2}\mo+\mo^{'}\mpq
\end{equation}
Moreover, if we form the difference:
$$A_{p-1}-\as_{p-1}=\frac{(p-1)!^2-1}{(p-1)!},$$
we see that the numerator of the right hand side is divisible by $p$. \\Hence, we have:
$$A_{p-1}-\as_{p-1}=(\as_{p-1}\;\mpt)(A_{p-1}^2-1)\mpq$$
Further, by using again Theorem $4$ with $2n=p-1$, we derive:
\begin{equation}A_{p-1}^2-1=p^2\,\mb_{p-1}^2-1-p^3\mc\mb(p-1)\mb_{p-1}+\frac{p^3}{3}\mbt(p-1)+\frac{5}{4}p^3\mb_{p-3}\mpq\end{equation}
And $\frac{1}{(p-1)!}\;\mpt$ is known from \cite{LEV3}, namely
$$\frac{1}{(p-1)!}=3(p\mb_{p-1}-1)-p\mb_{2(p-1)}-\frac{1}{2}p^2\mb_{p-1}^2\;\mpt$$
We thus get:
\begin{equation}\begin{split}
&A_{p-1}-\as_{p-1}=p^3(\mc\mb(p-1))_1-\frac{p^3}{3}\mbt(p-1)
-\frac{5}{4}p^3\mb_{p-3}\\&+\Big(-p\mb_{p-1}\big(p^2\mc\mb(p-1)\mpt\big)+p^2\mb_{p-1}^2-1\Big)\Big(\frac{1}{(p-1)!}\mpt\Big)\\
&\qquad\qquad\qquad\qquad\qquad\qquad\qquad\qquad\qquad\qquad\qquad\qquad\qquad\qquad\mpq\end{split}\end{equation}
with $p^2\mc\mb(p-1)$ modulo $p^3$ as provided in Result $11$ of the introduction. \\
By using Congruence $(84)$ for the first row of $(87)$ we finally derive:
\begin{equation}\begin{split}
A_{p-1}+\as_{p-1}&=p^2\big(\mc\mb(p-1)\big)_0-2p\mb_{p-1}\\&+\Big(p\mb_{p-1}\big(p^2\mc\mb(p-1)\mpt\big)-p^2\mb_{p-1}^2+1\Big)\Big(\frac{1}{(p-1)!}\,mod\,p^3\Big)\\
&\qquad\qquad\qquad\qquad\qquad\qquad\qquad\qquad\qquad\qquad\;\;\;\;\;\,\mpq
\end{split}\end{equation}
The final step consists of subtracting $(85)$ to $(88)$. It yields:
\begin{equation}\begin{split}
(p-1)!&=\frac{1}{3}\Bigg\lbrace 2p^2(2p\mb_{2(p-1)}-p^2\mb_{p-1}^2)_2-4p\mb_{p-1}-\frac{5}{4}p^3\mb_{p-3}-\frac{p^3}{2}\mo-\mo^{'}\\
%&-\Big(p^3\mb_{p-1}^3+p^2\mb_{p-1}^2-2(p\mb_{2(p-1)})(p\mb_{p-1})-1\Big)\Big(3(p\mb_{p-1}-1)-p\mb_{2(p-1)}-\frac{1}{2}p^2\mb_{p-1}^2\Big)\Bigg\rbrace\\
&-\Big(-p\mb_{p-1}\big(p^2\mc\mb(p-1)\mpt\big)+p^2\mb_{p-1}^2-1\Big)\Big(\frac{1}{(p-1)!}\,mod\,p^3\Big)\\
&\qquad\qquad\qquad\qquad\qquad\qquad\qquad\qquad\qquad\qquad\qquad\qquad\mpq
\end{split}
\end{equation}
Some additional and straightforward computations lead to the formula given in Theorem $3$.\\
Similarly to \cite{LEV3} (cf. Result $11$ of $\S\,2$), it is worth noting from $(84)$ that:
\begin{equation}\begin{split}
p^3\mbt(p-1)\!=3\big(\mc\mb(p-1)\big)_1p^3\!-6A_{p-1}\!+3\big(\mc\mb(p-1)\big)_0\!-\!6p\mb_{p-1}&\!\!-\frac{15}{4}p^3\mb_{p-3}\\&\mpq
\end{split}\end{equation}
with $A_{p-1}=(p-1)!$ modulo $p^4$ like provided in $(89)$, thus extending the results of Theorem $1$. \\

On a final note, it would be interesting to try and compute $(p-1)!$ mod $p^4$ by using the method of \cite{LEV1}, that is by using Hensel's algorithm in order to lift the $p$-adic integer roots of the polynomial $X^{p-1}+(p-1)!\in\mathbb{Z}_p[X]$ and draw a comparison with the result obtained here.\\

\textsc{Email address:} \textit{clairelevaillant@yahoo.fr}
\newpage
\begin{Large}Appendix A\end{Large}\\\\

$\mo=\mo_0+\mo_1+\mo_2+\mo_3+\mo_4+\mo_5+\mo_6\;\mpu$ with:
$$\begin{array}{cccc}
\mo_0&=&\bigg\lbrace 2\big(2p\mb_{2(p-1)}-p^2\mb_{p-1}^2\big)_2+(\mb_{p-1})_0^2+2(\mb_{p-1})_1+\big(2(\mb_{p-1})_0\big)_1\bigg\rbrace_1&\mpu\\
&&&\\
\mo_1&=&2(\mb_{p-1})_1(\mb_{p-1})_0+2(\mb_{p-1})_2-6\mb_{p-3}\big((\mb_{2})_1-(\mb_{p+1})_1\big)&\mpu\\
&&&\\
\mo_2&=&2(2p\mb_{2(p-1)}-p^2\mb_{p-1}^2)_2&\mpu\\
&&&\\
\mo_3&=&-\Big[2\Big[\frac{p+1}{2}(2p\mb_{2(p-1)}-p^2\mb_{p-1}^2)_2\Big]_0\Big]_1&\\&&&\\
&&-2\Big[\frac{p-1}{2}(2p\mb_{2(p-1)}-p^2\mb_{p-1}^2)_2\Big]_1-2\Big[2\Big[\frac{p-1}{2}(2p\mb_{2(p-1)}-p^2\mb_{p-1}^2)_2\Big]_0\Big]_1
&\mpu\\
&&&\\
\mo_4&=&-\Bigg\lbrace\Big[2\Big(1+\big(\mb_{2(p-1)}\big)_0\Big)\Big]_1
+1+2\Big(1+\big(\mb_{2(p-1)}\big)_0+\big(\mb_{2(p-1)}\big)_1\Big)\Bigg\rbrace_1&\\
&&-\Big[2\Big(1+\big(\mb_{2(p-1)}\big)_0\Big)\Big]_2-2\Big(1+\big(\mb_{2(p-1)}\big)_1+\big(\mb_{2(p-1)}\big)_2+2\big(1+\big(\mb_{2(p-1)}\big)_0\big)\Big)&\mpu\\
&&&\\
\mo_5&=&\big(3+2(\mb_{p-1})_0\big)_1+5+4(\mb_{p-1})_0+(\mb_{p-1})_0^2+2(\mb_{p-1})_1-\Big[p - 1 + \big(2 + 2(\mb_{p-1})_0\big)_0\Big]_2&\mpu\\
&&&\\
\mo_6&=&-\Bigg\lbrace(2p\mb_{2(p-1)}-p^2\mb_{p-1}^2)_2+\big(-3-2(\mb_{p-1})_0\big)_0+\Big[p - 1 + \big(2 + 2(\mb_{p-1})_0\big)_0\Big]_1\\
&&+\Big[\big[2(1+(\mb_{2(p-1)})_0)\big]_1+1+2\big(1+(\mb_{2(p-1)})_0+(\mb_{2(p-1)})_1\big)\Big]_0\Bigg\rbrace_1&\mpu
\end{array}$$

\begin{equation*}\begin{split}\mo^{'}=&\frac{p^2}{2}\Big(1+p\big((pB_{p-1})_1-1\big)\Big)(2p\mb_{2(p-1)}-p^2\mb_{p-1}^2)_2\\
&-4+6\,p\mb_{p-1}-8\,p\mb_{2(p-1)}+3\,p\mb_{3(p-1)}\\
%&+\frac{13p^3-12}{3}-\frac{19p^3-12}{2}p\mb_{p-1}-(8-14p^3)p\mb_{2(p-1)}-\frac{11p^3-6}{2}p\mb_{3(p-1)}\\
&+p^2\mb_{p-3}\Big(\frac{1}{4}-3(\mb_2)_0-p\big(\frac{1}{4}+3(\mb_{p+1})_1\big)\Big)\qquad\mpq\end{split}\end{equation*}
%&+\frac{p^3}{6}(\mb_{2p-4}-\mb_{p-3})_1-\frac{p^3}{30}(\mb_{2p-6}-\mb_{p-5})_1
%&-\frac{p^3}{3}(\mb_{2p-4}-\mb_{p-3})_1-\frac{p^3}{30}(\mb_{2p-6}-\mb_{p-5})_1\\
%&+\frac{p^2}{2}\mb_{p-3}\Big(1+p(1-6(\mb_{p+1})_1)\Big)-\frac{p^2}{30}\mb_{p-5}\Big(1+p(1+3(pB_{p-1})_1)\Big)

Let $$\widehat{\mo}:=\sum_{\begin{array}{l}0\leq i\leq 6\\i\neq 2\end{array}}\widehat{\mo_i}\mpu$$ with
$$\begin{array}{ccccc}\widehat{\mo_i}&=&\mo_i&\mpu &\forall i\neq 1,2\\
\widehat{\mo_1}&=&2(\mb_{p-1})_1(\mb_{p-1})_0+2(\mb_{p-1})_2&\mpu
\end{array}$$

\newpage

\begin{Large} Appendix B. Mathematica program. \end{Large}
\begin{verbatim}
wq[p_] := Mod[(1 + (p - 1)!)/p, p];                 (*Residue of the Wilson quotient*)
q[a_, p_] := Mod[(a^(p - 1) - 1)/p, p];             (*Residue of the Fermat quotient
                                                    in base a*)
dvB[b_] := BernoulliB[b]/b;                         (*Divided Bernoulli numbers*)
den[c_] := Denominator[c]; num[d_] := Numerator[d]; (*Abbreviations for numerators
                                                    and denominators*)

pr0[a_, p_] := Mod[Mod[num[a], p]*PowerMod[den[a], -1, p], p];
                                                    (*First residue in the p-adic
                                                    expansion of a p-adic integer*)

pr1[a_, p_] := pr0[(a-pr0[a,p])/p,p];               (*Second residue*)
pr2[a_, p_] := pr0[(a - pr0[a, p] - pr1[a, p]*p)/p^2, p];
                                                    (*Third residue*)
redt[a_,p_]:=Mod[Mod[num[a],p^3]*PowerMod[den[a],-1,p^3],p^3];
                                                    (*Reduction modulo p^3Z_p of a
                                                    p-adic integer*)
c[k_] := Sum[dvB[i]*dvB[k - i], {i, 2, k - 2}];     (*Convolution of divided Bernoulli
                                                    numbers of order k*)
ct[k_] := Sum[dvB[i]*c[k - i], {i, 2, k - 4}];      (*Triple convolution of divided
                                                    Bernoulli numbers of order k*)

dvBpz[p_] := pr0[dvB[p-1]-1/p,p];                   (*Second residue in the Hensel
                                                    of the p-adic number
                                                    the (p-1)-th dvB*)
dvBpu[p_] := pr0[(dvB[p - 1] - 1/p - dvBpz[p])/p, p];
                                                    (*Third residue*)
dvBpd[p_] := pr0[(dvB[p - 1] - 1/p - dvBpz[p] - p*dvBpu[p])/p^2, p];
                                                    (*Fourth residue*)
dvBdpz[p_] := pr0[dvB[2*(p - 1)] - PowerMod[2, -1, p]/p, p];
                                                    (*Second residue
                                                    2(p-1)-th dvB*)
dvBdpu[p_] := pr0[(dvB[2*(p - 1)] - PowerMod[2, -1, p]/p - dvBdpz[p])/p,p];
                                                    (*Third residue*)
dvBdpd[p_] := pr0[(dvB[2*(p - 1)] - PowerMod[2, -1, p]/p - dvBdpz[p]
- p*dvBdpu[p])/p^2, p];                             (*Fourth residue*)

otzh[p_] := pr1[2*pr2[2*p*dvB[2*(p - 1)] - p^2*dvB[p - 1]^2, p]+(*Other term OT_0 hat*)
   dvBpz[p]^2 + 2*dvBpu[p] + pr1[2*dvBpz[p], p], p];
otuh[p_] := pr0[2*dvBpu[p]*dvBpz[p] + 2*dvBpd[p], p];            (*Other term OT_1 hat*)
otth[p_] :=                                                      (*Other term OT_3 hat*)
  pr0[-2*pr1[(p - 1)/2*pr2[2*p*dvB[2*(p - 1)] - p^2*dvB[p - 1]^2, p],
     p] - pr1[
    2*pr0[(p + 1)/2*pr2[2*p*dvB[2*(p - 1)] - p^2*dvB[p - 1]^2, p], p],
     p] - 2*pr1[
     2*pr0[(p - 1)/2*pr2[2*p*dvB[2*(p - 1)] - p^2*dvB[p - 1]^2, p],
       p], p], p];
otqh[p_] :=                                                      (*Other term OT_4 hat*)
  pr0[-pr1[
     pr1[2*(1 + dvBdpz[p]), p] + 1 + 2*(1 + dvBdpz[p] + dvBdpu[p]),
     p] - pr2[2*(1 + dvBdpz[p]), p] -
   2*(1 + dvBdpu[p] + dvBdpd[p] + 2*(1 + dvBdpz[p])), p];
otch[p_] :=                                                      (*Other term OT_5 hat*)
  pr0[pr1[3 + 2*dvBpz[p], p] + 5 + 4*dvBpz[p] + dvBpz[p]^2 +
    2*dvBpu[p]-pr2[p - 1 + pr0[2 + 2*dvBdpz[p], p], p], p];
otsh[p_] :=                                                      (*Other term OT_6 hat*)
  pr0[-pr1[
     pr2[2*p*dvB[2*(p - 1)] - p^2*dvB[p - 1]^2, p] + pr0[- 3 - 2*dvBpz[p],p] + 
     pr1[p - 1 + pr0[2 + 2*dvBdpz[p], p], p] +
      pr0[pr1[2*(1 + dvBdpz[p]), p]+1+2*(1 + dvBdpz[p] + dvBdpu[p]),p], p],
   p];

oth[p_] := otzh[p] + otuh[p] + otth[p] + otqh[p] + otch[p] + otsh[p];
                                                                  (*Other term OT hat*)
fin[p_] :=                                                        (*RHS of Main Thm 3*)
  1/3*(-p^3*dvB[p-3]+p^2*PowerMod[2, -1, p^2]*
      pr2[2*p*dvB[2*(p - 1)] - p^2*dvB[p - 1]^2,
       p]*(3 - p*(1 + pr0[(p*BernoulliB[p - 1] + 1)/p,p])) -
     p^3*PowerMod[2, -1, p]*pr0[oth[p],p] + 4 - 10*p*dvB[p - 1] +
     8*p*dvB[2*(p - 1)] -3*p*dvB[3*(p - 1)] +
     (p*dvB[p - 1]*redt[2*p*dvB[2*(p-1)]-p^2*dvB[p-1]^2,p] - p^2*dvB[p - 1]^2 +1)
        *redt[3*(p*dvB[p - 1] - 1) - p*dvB[2*(p - 1)] -
        1/2*p^2*dvB[p - 1]^2,p]);

redfin[p_] :=
  Mod[Mod[num[fin[p]], p^4]*PowerMod[den[fin[p]], -1, p^4], p^4];(*Reduction mod p^4
                                                                 of RHS of Main Thm 3*)

  *Now testing the result by running the program on various random primes*

In:= Mod[18!, 19^4]                                              *Test on p=19*
Out= 93175

In:= redfin[19]
Out= 93175

In:= Mod[46!,47]                                                 *Test on p=47*
Out= 2266715

In:= redfin[47]
Out= 2266715

In:= Mod[60!, 61^4]                                             *Test on p=61*
Out= 6504002

In:= redfin[61]
Out= 6504002

In:= Mod[172!, 173^4]                                           *Test on p=173*
Out= 438178897

In:= redfin[173]
Out= 438178897

In:= Mod[520!,521^4]                                            *Test on p=521*
Out= 589386980

In:= redfin[521]
Out= 589386980

In:= Mod[876!,877^4]                                            *Test on p=877*
Out= 557572214137

In:= redfin[877]
Out= 557572214137

In:= Mod[1008!,1009^4]                                          *Test on p=1009*
Out= 709347287962

In:= redfin[1009]
Out= 709347287962

In:= NextPrime[10009]
Out= 10037

In:= Mod[10036!,10037^4]                                        *Test on p=10037*
Out= 3241073122386671

In:= redfin[10037]
Out= 3241073122386671

In:= Mod[120010!,120011^4]                                      *Test on p=120011*
Out= 143693568124824551692

In:= redfin[120011]
Out= 143693568124824551692
\end{verbatim}

\end{document}